\documentclass[3p,times]{amsart}

\usepackage{hyperref}
\usepackage{amsmath, amssymb}
\usepackage{amscd}
\usepackage{xcolor}
\usepackage{amsthm}
\usepackage[font=footnotesize]{caption}

\usepackage{enumerate} % for enumerate [(i)]/[(a)]/[(1)]

\newtheorem{theorem}{Theorem}[section] 
\newtheorem{lemma}[theorem]{Lemma}   
\newtheorem{corollary}[theorem]{Corollary}
\newtheorem{proposition}[theorem]{Proposition}
\newtheorem{definition}[theorem]{Definition}

\newtheorem{main-theorem}[theorem]{Theorem}
\newtheorem{assumption}[theorem]{Assumption}
\newtheorem*{problem*}{Problem}

\theoremstyle{definition}

\newtheorem*{question*}{Question}

\newtheorem{example}[theorem]{Example}
\newtheorem{remark}[theorem]{Remark}

\newcommand{\ul}{\underline}

\newcommand{\bN}{\mathbb{N}}

\newcommand{\bZ}{\mathbb{Z}}

\newcommand{\cO}{\mathcal{O}}
\newcommand{\cC}{\mathcal{C}}

\newcommand{\cT}{\mathcal{T}}

\newcommand{\cB}{\mathcal{B}}
\newcommand{\cQ}{\mathcal{Q}}

\newcommand{\cS}{\mathcal{S}}
\newcommand{\La}{\Lambda}
\newcommand{\ba}{\bar{\alpha}}

\newcommand{\ve}{\varepsilon}
\newcommand{\wt}{\widetilde}

\makeatletter
\newcommand{\vect}[1]{%
  \vbox{\m@th \ialign {##\crcr
  \vectfill\crcr\noalign{\kern-\p@ \nointerlineskip}
  $\hfil\displaystyle{#1}\hfil$\crcr}}}
\def\vectfill{%
  $\m@th\smash-\mkern-7mu%
  \cleaders\hbox{$\mkern-2mu\smash-\mkern-2mu$}\hfill
  \mkern-7mu\raisebox{-3.81pt}[\p@][\p@]{$\mathord\mathchar"017E$}$}

\newcommand{\amsvect}{%
  \mathpalette {\overarrow@\vectfill@}}
\def\vectfill@{\arrowfill@\relbar\relbar{\raisebox{-3.81pt}[\p@][\p@]{$\mathord\mathchar"017E$}}}

\newcommand{\amsvectb}{%
  \mathpalette {\overarrow@\vectfillb@}}
\newcommand{\vecbar}{%
  \scalebox{0.8}{$\relbar$}}
\def\vectfillb@{\arrowfill@\vecbar\vecbar{\raisebox{-4pt}[\p@][\p@]{$\mathord\mathchar"017E$}}}
\makeatother
\usepackage{etex} % need because of using both tikz and xymatrix below

\usepackage{etex} % needed because of using both tikz and xymatrix below
\usepackage{dsfont}         % for \mathds{1}
\usepackage[h]{esvect}

\usepackage[matrix,arrow,curve,frame,arc]{xy}

\usepackage{pgf}
\usepackage{tikz}

\usepackage{color}

\usetikzlibrary{arrows,shapes,automata,backgrounds,petri,calc}

\newcommand{\tikzAngleOfLine}{\tikz@AngleOfLine}
\def\tikz@AngleOfLine(#1)(#2)#3{%
\pgfmathanglebetweenpoints{%
\pgfpointanchor{#1}{center}}{%
\pgfpointanchor{#2}{center}}
\pgfmathsetmacro{#3}{\pgfmathresult}%
}

\begin{document}

\baselineskip=14pt

\title{Hybrid algebras}

{\def\thefootnote{}
\footnote{This research was supported by the program
``Research in Pairs'' by the 
Mathematisches Forschungsinstitut Oberwolfach in 2018, and also by the
	Faculty of Mathematics and Computer Science of the Nicolaus Copernicus University in Toru\'{n}. Work on this paper
	was in progress when in October 2020, sadly, Andrzej passed away.}
}

\author[K. Edrmann]{Karin Erdmann}
\address[Karin Erdmann]{Mathematical Institute,
   Oxford University,
   ROQ, Oxford OX2 6GG,
   United Kingdom}
\email{erdmann@maths.ox.ac.uk}

\author[A. Skowro\'nski]{Andrzej Skowro\'nski}
\address[Andrzej Skowro\'nski]{Faculty of Mathematics and Computer Science,
   Nicolaus Copernicus University,
   Chopina~12/18,
   87-100 Toru\'n,
   Poland}

\begin{abstract}
We introduce a new class of symmetric algebras, which we call hybrid algebras. This class   contains on
one extreme Brauer graph algebras, and
on the other extreme  general weighted surface algebras.
We show that hybrid algebras are precisely the blocks of idempotent algebras
of weighted surface algebras, up to socle deformations. 
More generally, for tame symmetric  algebras whose Gabriel quiver is 2-regular, we show that the tree class of an arbitrary Auslander-Reiten component is Dynkin or
Euclidean or one of the infinite trees $A_{\infty}, A_{\infty}^{\infty}$ or $D_{\infty}$.
\bigskip

 

\end{abstract}

%\linenumbers

\maketitle

\bigskip

\section{Introduction}\label{sec:Intro}

We are interested in the representation theory of tame self-injective algebras.
In this paper,  all algebras are finite-dimensional basic associative and  
indecomposable  $K$-algebras  
over an algebraically closed field $K$ of arbitrary characteristic.

In  the  modular representation theory of finite groups
representation-infinite tame blocks occur only
over fields of characteristic 2,  and their defect groups are dihedral,
semidihedral, or  (generalized) quaternion 2-groups.
Such blocks were studied  in a more general setting:  this led to  algebras of dihedral, 
semidihedral and quaternion type,  over algebraically closed fields of
arbitrary characteristic, which were introduced and investigated in \cite{E1}. 
These algebras are quite  restrictive, for example the number
of simple modules can be at most $3$, and one would like to know
how these fit into a wider context.

Recently cluster theory has led to new directions. 
Inspired by  this, 
we study in \cite{WSA}, \cite{WSA-GV}, \cite{WSA-corr} and \cite{WSA-SD} 
a  class of symmetric algebras defined in terms of surface
triangulations, which we call \emph{weighted surface algebras}.  They are periodic as algebras
of period 4 (with a few exceptions).
All  but one  of the algebras of quaternion type occur in this setting. Furthermore, most algebras of dihedral type, and of semidihedral type 
occur naturally as degenerations of these weighted surface algebras. As well, Brauer graph algebras, which includes
blocks of finite type, appear. This places  blocks of finite or tame representation
type into a  much wider context, which also connects with other parts of mathematics.

In this paper, we present a unified approach. We 
introduce  a new class of algebras, which we call {\it hybrid algebras}.
At one extreme it contains all Brauer graph algebras, and  at the other  extreme it contains all 
weighted surface algebras, which are almost all
 periodic as algebras, of period four (see \cite{WSA} and \cite{WSA-GV}). 
Furthermore, the class  contains  many other symmetric algebras of tame or finite representation type. 
In particular it contains all blocks of group algebras, or type $A$ Hecke algebras, of tame or finite type, up to Morita equivalence.

Consider tame symmetric algebras more generally. One observes that being tame is a strong restriction on the Gabriel quiver of the algebra. 
At any given vertex there are not too many arrows starting or ending, and also not too few, avoiding finite type. The situation when one can expect
larger classes of algebras occurs when the Gabriel quiver is 2-regular. We ask whether all tame symmetric algebra with a 2-regular Gabriel quiver are hybrid algebras,
up to some small list of exceptions, and up to derived equivalence.  Our result on the tree class of stable AR components holds for any tame symmetric algebra
with 2-regular Gabriel quiver, and could be thought of as some evidence.

A motivation is that  various basic 
tame, or finite type, 
symmetric  algebras studied in recent years
 have a unified description, of the form $\La = KQ/I$  with $(Q, I)$ satisfying certain combinatorial restrictions. Namely, 
 the quiver $Q$ is 2-regular, that is,
there are two arrows starting and two arrows ending at each vertex. 
Here $I$ may contain arrows of $Q$,  so that the Gabriel quiver can be seen 
as a subset of $Q$. The fact that $Q$ is 2-regular, gives rise to symmetry.  
There is  a permutation $f$ of the arrows such that
 $t(\alpha) = s(f(\alpha))$ for each arrow $\alpha$.
This  determines uniquely a different permutation $g$ where
$t(\alpha) = s(g(\alpha))$ but $f(\alpha)\neq g(\alpha)$.
Such permutations have been studied for Brauer graph algebras: the permutation $g$ describes the cyclic order in  the Brauer graph, and the permutation $f$ has
been called the 'Green walk'. Here we will see that these permutations
$f$ and $g$ exist more generally. 

The permutation $f$ encodes   minimal relations, and 
the permutation $g$ describes, roughly speaking, a basis for
the indecomposable projective modules.
Consider $e_i\La$, and let $\alpha, \ba$ be the arrows starting at $i$.
Then $e_i\La$ has a basis consisting of monomials along the $g$-cycles
of $\alpha$ and of $\ba$, and the socle of $e_i\La$ is spanned
by  $B_{\alpha}$ (or $B_{\ba}$), where $B_{\alpha}$ is  the longest monomial
 starting with $\alpha$ which is non-zero in $\La$. Let also $A_{\alpha}$ be the submonomial of $B_{\alpha}$
such that $B_{\alpha} = A_{\alpha}\gamma$ where $\gamma$ is the arrow with $g(\gamma)=\alpha$. 

\medskip
For  each  arrow $\alpha$ there is a minimal relation
 determined by $f$, either  'biserial', or 'quaternion':\\ 
(B) \ $\alpha f(\alpha) \in I$, or\\
(Q) \ $\alpha f(\alpha) - c_{\ba}A_{\ba} \in I$ \\
(where the $c_{\ba}$ are non-zero scalars constant on $g$-cycles). \ 
With these data,  
together with suitable zero relations, and up to socle deformations,
the following hold.

 \ The algebra
$\La$ is a Brauer graph algebra if all minimal relations are biserial.
If $f^3=1$ and all minimal relations are quaternion, then the algebra $A$ is a weighted surface algebra (as in \cite{WSA, WSA-GV, WSA-corr}).
When   $f^3=1$,  and some but not all minimal relations are biserial, we get algebras generalizing  
algebras of semidihedral type, as in \cite{E1} (see also \cite{Haj}).
As well algebras of finite type can occur naturally (which we also call tame in this context).

The known structure of tame local symmetric algebras should be further
motivation.
As one can find in \cite{E1}, section III, up 
to socle deformations, only relations of
the form (B) and (Q) occur. This suggests that 'generally' it should be sufficient to incorporate these types of relations.
Cycles of $f$ of length 3 (or 1) play a special role in the algebras studied in \cite{E1}. A relation (Q) only occurs if $\alpha$ belongs such a cycle. 
Namely we have $A_{\ba} g^{-1}(\ba) = B_{\ba}$ and $g^{-1}(\ba) = f^{-1}(\alpha)$ therefore
$\alpha f(\alpha)f^{-1}(\alpha)$ is a cyclic path, so the arrow $\alpha$ occurs in some triangle. 
\bigskip

We call the  set of arrows in an $f$-cycle of length $3$ or $1$ a {\it triangle}. 
Describing a {\it hybrid algebra $H$} in broad terms, we fix a set $\cT$ of triangles in $\cQ$.
Then $H= H_{\cT} = KQ/I$ where\\
(i) an arrow $\alpha \in \cT$ satisfies the quaternion relation, and\\
(ii) an arrow $\alpha \not\in \cT$ satisfies the biserial relation. \\
In addition there are zero relations. 

We start with a hybrid algebra where the quiver $Q$ for the definition is 
the Gabriel quiver, this is introduced and studied in Section 2. We  call the
algebras {\it regular} hybrid algebras. This is extended in  Section 3.
Our first main result is the following.

\begin{theorem}\label{thm:main1}  
(i) Assume $\La$ is a weighted surface algebra and $e$ is an idempotent of $\La$, then every block  component
	of $e\La e$ is a hybrid algebra (up to socle equivalence).

(ii) Assume $H$ is a hybrid algebra. Then there is a weighted surface
algebra $\La$ and an idempotent $e$ of $\La$ such that $H$ is isomorphic to a block component
of $e\La e$.\\
\end{theorem}

\medskip

The second part of this theorem  generalises \cite{BGA} where  we prove
that every Brauer graph algebra occurs
as an idempotent algebra of a weighted surface algebra.
For the second part, given a hybrid algebra $H$,  to construct the weighted surface algebra
$\La$ with $H$ as a component of $e\La e$, we use the $*$ construction
introduced in \cite{BGA}. 

Idempotent algebras of weighted surface
algebras  include many local algebra,  therefore our definition of
hybrid algebras must included these. In our general construction of
weighted surface algebras \cite{WSA-GV}, we have allowed virtual arrows, with
the benefit of essentially enlarging the class of algebras. The
price to pay is that zero relations have to be treated with care
(see  \cite{WSA-corr}), and naturally this
is also the case for hybrid algebras. In particular we need to
exclude a few small algebras
(see  Assumption \ref{ass:3.4}).

All local symmetric algebras 
of tame or finite type, and 
almost all algebras of dihedral, semidihedral or quaternion
type as in \cite{E1} are hybrid algebras.
There is one family of algebras of quaternion type which are not  hybrid algebras, but
are derived equivalent to  algebras of quaternion type (algebras $Q(3\cC)^{k, s}$, see \cite{H2}).

Hybrid algebras 
place blocks into a wider context; 
 in \cite{AGQT}  we define 
algebras of generalized
quaternion type, as  tame symmetric
algebras with periodic module categories,  that is,  generalizing quaternion
blocks, and show
that the ones with  2-regular Gabriel quiver are almost all
weighted surface algebras. 
As well in \cite{AGDT} we define algebras
of generalized dihedral type, in terms of homological properties
generalizing dihedral blocks, and show that almost all  are
the biserial weighted surface
algebras as in \cite{WSA}.
One would like a similar homological description of the hybrid algebras
which generalize semidihedral blocks.  

In order to understand the representation theory for  all  these algebras, the structure of the stable Auslander-Reiten quiver is 
essential. Our second main result is more general, it describes its graph structure for arbitrary
tame symmetric algebras with 2-regular Gabriel quiver:

\medskip

\begin{theorem} \label{thm:main2} Assume $\La$ is a tame symmetric algebra with a 2-regular  Gabriel quiver. Then the tree classes 
of stable Auslander-Reiten components of $\La$ are one of the infinite trees $A_{\infty}, A_{\infty}^{\infty}$ or $D_{\infty}$, or Euclidean or Dynkin.
\end{theorem}

It would be interesting to know whether a component with tree class $A_{\infty}$ of a tame symmetric
algebra  is necessarily a tube.

\medskip

We describe the organisation of the paper. 
In Section 2, we present and study a simplified version of hybrid algebras, which we call regular. For such an algebra,  $Q$ is the Gabriel quiver. In this case
we   prove 
a weaker version of Theorem 1.1, which will show how virtual arrows 
occur. 

In Section 3  we give the general definition, 
and discuss exceptions for the zero relations.
The details for consistency and bases are
refinements of results in Section 2 and are therefore
only given in an appendix. Originally we had incorporated socle deformations
into the general definition of a hybrid algebras. This is not done here, as it has caused further technical work. Note however that socle deformations can occur but are easy to identify.

In Section 4 we 
discuss algebras with few simple modules and small multiplicities. In Section 5 we prove Theorem \ref{thm:main1}, extending the
version in Section 2.   Section 6
is valid more generally, for arbitrary tame symmetric algebras
with 2-regular Gabriel quiver. The main result is Theorem
\ref{thm:main2} on stable Auslander-Reiten components. In the case of
hybrid algebras, we identify  components containing simple modules, and see
in particular that
the infinite trees in the list all   occur.

For further background and motivation, we refer to  \cite{ASS, Benson}, and to the introductions of \cite{WSA, WSA-GV},  or \cite{BGA}.

\bigskip

\section{Preliminaries and regular hybrid algebras}\label{sec:def}

\medskip

\subsection{The setup}\label{ssec:2.1}

Recall that a quiver is a quadruple $Q = (Q_0, Q_1, s, t)$ 
where $Q_0$ is a finite set of vertices, 
$Q_1$ is a finite set of arrows, and 
where $s, t$ are maps $Q_1\to Q_0$ associating 
to each arrow $\alpha\in Q_1$ its source $s(\alpha)$  
and its target $t(\alpha)$. 
We say that
$\alpha$ starts at $s(\alpha)$ and ends at $t(\alpha)$. 
We assume throughout that any quiver is connected. 
The quiver $Q$ is 2-regular if at each vertex, two arrows start and
two arrows end.

Denote by $KQ$ the path algebra of $Q$ over $K$. 
The underlying space has basis given by the set of all paths in $Q$, in particular 
for each vertex $i$, let $\ve_i$ be the path of length zero at $i$ in $KQ$.  
We will consider algebras of the form $\La = KQ/I$ for some ideal $I$ of $KQ$. Let $e_i= \ve_i +I$, then
the $e_i$ are pairwise orthogonal idempotents, and their sum is the identity of $\La$. 
We assume that the ideal $I$  contains all paths of length $\geq N$ for  some $N\geq 2$, so that
the algebra is finite-dimensional and basic. 
The Gabriel quiver $Q_{\La}$  of $\La$ has by definition the same vertices as $Q$ and its arrows
are in bijection with a basis for $J/J^2$ where $J$ is the radical of $\La$. Usually, $Q_{\La}$ 
can be taken as a subquiver of $Q$.

\subsection{Notation}\label{ssec:2.2} 
Recall that a biserial quiver  is a pair $(Q, f)$ 
where $Q$ is a 2-regular quiver, and $f$ is a permutation of the arrows 
such that 
for each arrow $\alpha$ we have $s(f(\alpha)) = t(\alpha)$.
This was defined in \cite{BGA}, but  here we also
allow the quiver $Q$ with only one vertex.  Moreover, we 
have an involution $(-)$ on the arrows, taking $\ba$ to be the arrow $\neq \alpha$ with the same starting vertex.
Given $f$, this uniquely determines the permutation $g$ on arrows, defined by $g(\alpha) = \overline{f(\alpha)}$.

Let $\mathcal{O}$ be the set of $g$-orbits on $Q_1$.
We fix a weight function (or multiplicity function), that is, a function  $m_{\bullet}: \cO(g)\to \bN$, 
and we fix a parameter function, that is, a function  $c_{\bullet}: \cO(g) \to K^*$. 
Moreover, $n_{\alpha}$ is the size of the $g$-orbit of   $\alpha\in Q_1$.

For an arrow $\alpha$ of $Q$, let 
$B_{\alpha}$ be the monomial along the $g$-cycle of $\alpha$ which starts with $\alpha$, 
of length $m_{\alpha}n_{\alpha}$, and let
$A_{\alpha}$ be the  submonomial of $B_{\alpha}$ starting with $\alpha$ of length $m_{\alpha}n_{\alpha}-1$, so that
$B_{\alpha} = A_{\alpha}g^{-1}(\alpha)$. 

For a path $p$ in $KQ$ we write $|p|$ for the length of $p$.
We will sometimes
write $p\equiv q$ if  $p$ and $q$ are  paths  in $KQ$ such that $p= \lambda q$ in some algebra $KQ/I$ for $0\neq \lambda\in K$.

\bigskip

\subsection{Regular hybrid algebras}\label{ssec:2.3}

The arrows in   $f$-orbits of length 3 or 1 play a special role,  we refer to these as {\it triangles}.
Note that any set of triangles is invariant under the permutation $f$.
The regular hybrid algebra is defined so that it has $Q$ as its Gabriel quiver, this is ensured by the following:

(*) \ We assume $m_{\alpha}n_{\alpha}\geq 2$ for any arrow $\alpha$, and  $m_{\alpha}n_{\alpha}\geq 3$ if $\ba  \in  \cT$.

\medskip

\begin{definition}\normalfont \label{def:2.1}  Let $(Q, f)$ be a biserial quiver with  the data  $m_{\bullet}, c_{\bullet}$ as in \ref{ssec:2.2}, 
and let  $\cT$ be a set of distinguished triangles. 
The regular hybrid algebra  $H=H_{\cT} = H_{\cT}(Q, f, m_{\bullet}, c_{\bullet} )$ associated to
	$\cT$, with assumption (*), is the algebra 
$H= KQ/I$ where $I$ is generated by the following elements:

	(1) \ 
	$\alpha f(\alpha) - c_{\bar{\alpha}}A_{\bar{\alpha}}$\ for  $\alpha \in \cT$ \ \ and  $ \alpha f(\alpha)$ for  $\alpha\not\in\cT$.

	(2) \ $\alpha f(\alpha) g(f(\alpha))$  
	and 
	\  $\alpha g(\alpha)f(g(\alpha))$ for any arrow $\alpha$ of $Q$.
			
	(3) \ $c_{\alpha}B_{\alpha} - c_{\bar{\alpha}}B_{\bar{\alpha}}$ for any arrow $\alpha$ of $Q$. 

\end{definition}

Let $i$ be a vertex and $\alpha, \ba$ the arrows starting
at $i$.
We say that  $i$ is {\it biserial} if  $\alpha$ and $\ba$ are both not in $\cT$.
We call the vertex $i$ a  {\it quaternion} vertex if $\alpha$ and $\ba$ are both in $\cT$.
Otherwise, we say that $i$ is {\it hybrid}.

\bigskip
The conditions (*) 
imply that arrows are not contained in $I$, so that  $Q$ is the Gabriel quiver of $H$. 
If $\cT = \emptyset$, then the algebra
$H$ is special biserial and symmetric, that is,  a Brauer graph algebra (BGA). At the
other extreme, if $\cT = Q_1$ then $H$ is a weighted surface algebra (WSA), as defined in \cite{WSA},  if $Q$ has at least three vertices, 
or it occurs amongst the algebras of quaternion type in \cite{E1}.

\begin{example}\label{ex:2.2} Consider the quiver
\[
%  \xymatrix@R=2pc@C=1.5pc{
%  \xymatrix@R=3.5pc@C=1.8pc{
  \xymatrix@R=3.pc@C=1.2pc{
%  \xymatrix@C=.8pc{
    &&& i
    \ar[ld]^{\tau}
    \ar[rrrd]^{\alpha}
    \\
    k
    \ar[rrru]^{\beta}
    \ar@<-.5ex>[rr]_(.6){\omega}
    && y
    \ar@<-.5ex>[ll]_(.4){\xi}
    \ar[rd]^{\rho}
    &&&& j \ar@(ru, dr)^{\eta}
    \ar[llld]^{\delta}  
   \\
    &&& x
    \ar[lllu]^{\gamma}
    \ar@/_3ex/[uu]^{\sigma}
  }
\]
As the permutation $f$, we take
	$$f = (\alpha \ \delta \  \sigma)(\rho \ \gamma \ \omega) (\xi \ \beta \ \tau)(\eta)$$
Then 
$$g = (\alpha \ \eta \ \delta \ \gamma \ \beta)(\tau \ \rho \  \sigma ) (\xi \ \omega)$$
	We take  $m_{\alpha}=1=m_{\tau}$ and $m_{\xi} =2$ and $c_{\alpha}=c, c_{\ba}=d$ and $c_{\xi}=1$.

The permutation $f$ has four cycles, each of size $1$ or $3$, so there are
several choices for the set $\cT$ of distinguished triangles.\\ 
	(a) \ If $\cT = Q_1$ then the algebra $H_{\cT}$ is a weighted surface
	algebra, as in \cite{WSA}.\\
	(b) \ If $\cT = \emptyset$ then the algebra is  special biserial and symmetric, hence a Brauer graph algebra. \\
	(c) \ An example for an intermediate choice of $\cT$ might be $\cT = \{\alpha, \delta,  \sigma, \eta\}$.
Then the relations for the paths of length $2$ between arrows of $\cT$ are
	$$\alpha \delta = dA_{\tau}, \ \ \delta\sigma = cA_{\eta}, \ \  \sigma\alpha = c A_{\gamma}, \ \eta^2 = cA_{\delta};
$$
and products of paths of length two along each other $f$-cycle are zero in 
$H_{\cT}$. In this case, vertices $i$ and $x$ are hybrid, vertex $j$ is quaternion, and vertices $k, y$ are biserial.
\end{example}

\begin{lemma} \label{lem:2.3} The conditions (1) to (3) in Definition
	\ref{def:2.1} are  consistent. In particular 
	$B_{\alpha}$ is non-zero on $H$.
\end{lemma}

\begin{proof}   We  show that
the condition for $\alpha f(\alpha)$  from (1)  and  the conditions 
for $g^{-1}(\alpha)\alpha f(\alpha)$ and $\alpha f(\alpha)g(f(\alpha))$ from (2) agree.
This is clear when $\alpha\not\in \cT$ since then condition (1) requires
$\alpha f(\alpha)=0$ in $H$. 

Assume now that $\alpha\not\in \cT$, then 
we substitute $\alpha f(\alpha) = c_{\ba}A_{\ba}$. We should have that  
$g^{-1}(\alpha) A_{\ba} = 0$ in $H$. 
By the definition of the permutations, we have  $g^{-1}(\alpha) = f^{-1}(\ba)$, and 
by the  assumption  (*), the monomial  $A_{\ba}$ has length at least $2$ and therefore
$f^{-1}(\ba)A_{\ba} = f^{-1}(\ba)\ba g(\ba) p$ for some monomial $p\in KQ$ of length $\geq 0$. 
Now condition (2) gives that this is zero in $H$. 
Similarly,  $A_{\ba}g(f(\alpha)) = q g^{-1}(\beta) \beta f(\beta)$ where $\beta = g^{-2}(\ba)$ is the last arrow of $A_{\ba}$ and $q\in KQ$ a monomial
of length $\geq 0$, and this is zero in $H$ by condition (2).
Similarly one verifies that conditions (1) and (3) agree.
\end{proof}

\bigskip

\begin{lemma} \label{lem:2.4} For each vertex $i$ and arrow $\alpha$ starting at $i$, we have
$B_{\alpha}J=0$ and $JB_{\alpha}=0$ where $J$ is the radical of $H$. In particular $B_{\alpha}\neq 0$ belongs to the socle of
$e_iH$.
\end{lemma}

\begin{proof}  \ We have 
	$B_{\alpha}\alpha = \alpha B_{g(\alpha)} \equiv \alpha B_{\bar{g(\alpha)}}  = \alpha B_{f(\alpha)} 
	=  \alpha f(\alpha) g(f(\alpha)) p$ 
where $p$ is some monomial of length $\geq 0$ and this is zero by condition (2). Then we have as well that
$B_{\alpha}\ba \equiv B_{\ba}\ba = 0$. 
\end{proof}

We write $(B_{\alpha})_j$ for the initial submonomial
$\alpha g(\alpha)\ldots g^{j-1}(\alpha)$ of $B_{\alpha}$ of length $j$.

\bigskip

\begin{lemma} \label{lem:2.5} Let $\alpha \in Q_1$, and let 
	$\cB_{\alpha}:= \{ (B_{\alpha})_j \mid 1 \leq j\leq |B_{\alpha}|\}$ be
	the set of all initial submonomials of $B_{\alpha}$.
\begin{itemize}
\item [(a)]
 The set $\cB_{\alpha}$ is 
	 linearly independent in $H$. 
	\item[(b)]
Assume that  $\alpha, \ba$ are both in $\cT$,  then 
	$\cB_{\alpha} \cup A_{\ba}$ also is 
	linearly independent.
	\end{itemize}
\end{lemma}

\bigskip

\begin{proof}  \ (a) \ Let
$$\sum_{j=1}^{|B_{\alpha}|} a_j(B_{\alpha})_j = 0 \ \ (a_j\in K).
$$
	Premultiplying  with $A_{g^{-1}(\alpha)}$ gives
$0= a_1A_{g^{-1}(\alpha)} \alpha = a_1B_{g^{-1}(\alpha)}$ and hence $a_1=0$. 
Suppose we have $a_1=\ldots = a_{r-1}=0$. We premultiply with 
the submonomial $q$ of $B_{\alpha}$ such that $q (B_{\alpha})_r$ is equal to $B_{\gamma}$ for the appropriate $\gamma$. This annihilates all terms except one, leaving only $a_rB_{\gamma}=0$ and so $a_r=0$. 

(b) \ Let $\sum_{j=1}^{|B_{\alpha}|} a_j(B_{\alpha})_j + b A_{\ba}=0$ with $a_j$ and $b$ in $K$.  
We premultiply with $f^{-1}(\ba) = g^{-1}(\alpha)$. 
By condition (2) of  Definition \ref{def:2.1}, using also that $|A_{\ba}|\geq 2$ we get
$f^{-1}(\ba)A_{\ba} =0$, and this leaves 
$$\sum_{j=1}^{|B_{\alpha}|}a_j(B_{g^{-1}(\alpha)})_{j+1}=0. 
$$
Hence $a_1= \ldots = a_{|B_{\alpha}|-1}=0$, by (a), and we are
	left with $a_{|B_{\alpha}|}B_{\alpha} + bA_{\ba} = 0$. 
	Using that $B_{\alpha} \equiv B_{\ba}$, we have linear
	combination of two initial submonomials of $B_{\ba}$, and by part  (a)
	(applied to $\ba$),  the coefficients are zero.
\end{proof}

\bigskip

\begin{lemma}\label{lem:2.6}  The module $e_iH$ has basis
	$\{ e_i\} \cup \cB_{\alpha} \cup \cB_{\ba} \setminus \{B_{\ba}\}$. Hence $\dim e_iH = m_{\alpha}n_{\alpha} + m_{\ba}n_{\ba}$.
\end{lemma}

\begin{proof} \ Suppose we have
$$\sum_{j=1}^{|B_{\alpha}|} a_j(B_{\alpha})_j + \sum_{t=1}^{|B_{\ba}|-1}\bar{a}_t(B_{\ba})_t= 0. \ \leqno{(*)}
$$

(a) Assume first that (say) $\ba$ is not in $\cT$. We premultiply (*) with 
$f^{-1}(\ba)$, this annihilates the second sum. Recall $f^{-1}(\ba) = g^{-1}(\alpha)$, therefore the first sum becomes 
$$0 = \sum a_j(B_{g^{-1}(\alpha)})_{j+1},
$$ and by  Lemma \ref{lem:2.5},  $a_j=0$ for all $j< |B_{\alpha}|$. 
Then (*) becomes
$$0 = a_{|B_{\alpha}|}B_{\alpha} + \sum_{t=1}^{|B_{\ba}|-1} \bar{a}_t(B_{\ba})_t= 0.
$$
Since $B_{\alpha} \equiv B_{\ba}$ we can again apply Lemma \ref{lem:2.5} 
and deduce that all coefficients are zero.

(b) Assume  $\alpha, \ba$ are both in $\cT$. 
We premultiply with $\gamma = f^{-1}(\ba)$. We have $\gamma\ba = c_{\bar{\gamma}}A_{\bar{\gamma}}$ but $\gamma\ba g(\ba)=0$ and there is only one
non-zero term from the second sum, namely a multiple of $A_{\bar{\gamma}}$. The first sum is a linear combination of elements $(B_{\gamma})_j$ since $\gamma = g^{-1}(\alpha)$.
We apply  part (b) of Lemma \ref{lem:2.5}   and deduce that all scalar coefficients are zero.
\end{proof}

\bigskip

\subsection{Idempotent algebras of WSA's}\label{ssec:2.4}

In  \cite{WSA} we have
studied weighted surface algebras  whose 
Gabriel quiver  is 2-regular (with at least three vertices). 
One may ask whether 
  an idempotent algebra of such a WSA is a regular 
 hybrid algebra. 
We will investigate  this, and determine when exactly this
is the case, and at the same time it will illustrate why we should
allow virtual arrows  for  general hybrid
algebras. Examples can be found in 2.8 below.

\bigskip

  \begin{proposition} \label{prop:2.7} Assume $\La$ is a WSA with a 2-regular Gabriel quiver.
	  Let $\Gamma$ be a subset of $Q_0$ and $e=\sum_{i\in \Gamma} e_i$, and
	  let $R = e\La e$.\\
	  (i) The idempotent algebra $R$ satisfies conditions (1) to (3) of Definition \ref{def:2.1}.\\
	  (ii)  $R$ satisfies the multiplicity condition (*)  unless for some  $i\in \Gamma$ and $\alpha$ starting at $i$ we have
	  	   \begin{itemize}
	   \item [(*1)]  $m_{\alpha}=1$ and 
	  the $g$-cycle of $\alpha$ intersects $\Gamma$ only in $i$ (with no repetition);  or  
\item[(*2)]  $m_{\alpha} = 1$ and $\wt{n}_{\wt{\alpha}} = 2$, and $\Gamma$  contains both $s(g^{-1}(\alpha))$ and $t(\ba)$.
\end{itemize}
  \end{proposition}

  \bigskip

 \begin{proof} 
\ \ Let $\La$ be a WSA with 2-regular Gabriel quiver, that is it has a presentation
$\La = KQ/I$ 
of a (regular) hybrid algebra such that $\cT = Q_1$. 
In  particular we have then $m_{\alpha}n_{\alpha}\geq 3$ for all $\alpha$.
The only additional assumption in \cite{WSA} is that the quiver has 
at least three vertices  (see the text following \cite[Theorem 1.4]{WSA}).
Take a subset $\Gamma$ of $Q_0$, and let $e=\sum_{i\in \Gamma} e_i$ and $R:= e\La e$.

(i) \ We compute the basic algebra for $R$. Let $\wt{Q}$ be the quiver  with 
 vertices corresponding to  the primitive idempotents of $R$, that is  the $e_i (= ee_ie)$ with $i\in \Gamma$. 
 For $\alpha \in Q_1$ and $s(\alpha) = i\in \Gamma$, let 
$\wt{\alpha}$ be  the shortest path in $Q$ along the $g$-cycle of $\alpha$, starting with $\alpha$, and ending at a vertex in $\Gamma$.
We define $\wt{Q}$ by taking  the set 
$$\wt{Q}_1 = \{ \wt{\alpha} \in KQ\mid \alpha \in Q_1, \ \alpha =e_i\alpha
\mbox{ for } i\in \Gamma\}$$
as its set of arrows.   The set  $\wt{Q}$ is a   generating set for the radical of $R$, and hence 
we have a
surjective algebra map
$\psi: K\wt{Q} \to R$, and $R\cong K\wt{Q}/\wt{I}$ where $\wt{I}$ is
the kernel of $\psi$.

(a) \ We observe that
the quiver $\wt{Q}$ is 2-regular: We have two arrows starting at each vertex, and also two arrows ending at each vertex
(write $B_{\alpha}, B_{\ba}$ as a product of elements in $\wt{Q}_1$, then $B_{\alpha}$ and $B_{\ba}$ end with distinct arrows of $\wt{Q}$).

\ \ We define a permutation $\wt{f}$. 
Let $\wt{\alpha} = \alpha g(\alpha)\ldots g^{p}(\alpha)$ and $\beta := f(g^p(\alpha))$, then 
$$\wt{f}(\wt{\alpha}):= \wt{\beta}.
$$
With this, each connected component of $(\wt{Q}, \wt{f})$ is a biserial quiver.
Furthermore,  the permutation $\wt{g}$ is obtained
from the cycles of $g$ in $Q$, by factorizing them at each vertex in $\Gamma$.
In particular if $\wt{n}_{\wt{\alpha}}$ is the length of
the cycle of $\wt{\alpha}$, then $1\leq \wt{n}_{\wt{\alpha}} \leq n_{\alpha}$.
The multiplicity function $\wt{m}$ for $\wt{Q}$ 
must be taken as   $\wt{m}_{\wt{\alpha}} = m_{\alpha}$, and
the parameter function $\wt{c}$ is taken as  $\wt{c}_{\wt{\alpha}} = c_{\alpha}$ for each arrow $\wt{\alpha}$. 
Note that we may view the path algebra $K\wt{Q}$ as a subspace of $KQ$, and if so then $B_{\wt{\alpha}}$ is equal to $B_{\alpha}$. 

(b) \ There is a canonical set $\wt{\cT}$ of distinguished triangles of $\wt{Q}$. Let
 $$\wt{\cT}:= \{ \wt{\alpha}\mid  \alpha = \wt{\alpha} \mbox{ and } \wt{f(\alpha)} = f(\alpha)\}
 $$
 Note that if $\alpha =\wt{\alpha}$ and also $f(\alpha) = \wt{f(\alpha)}$ then 
 both $s(\alpha)$ and $t(f(\alpha))$ are  in $\Gamma$, and hence
 $f^2(\alpha) = \wt{f^2(\alpha)}$. Therefore  $\wt{\cT}$ is  closed under under the permutation $\wt{f}$.
 Furthermore, the arrows in $\wt{\cT}$ satisfy the  relations (1) of Definition \ref{def:2.1}.

(c) \  We show now that  for $\wt{\alpha}\not\in  \wt{\cT}$ we have
$\wt{\alpha}\wt{f}(\wt{\alpha}) = 0$. 
With the notation as in (a) we have
$$\wt{\alpha}\wt{f}(\wt{\alpha}) = \alpha g(\alpha)\ldots g^p(\alpha)f(g^p(\alpha))q
\leqno{(*)}$$
for some monomial $q\in KQ$. If $p\geq 1$ this is zero in $\La$, by condition (2) of Definition \ref{def:2.1}. 
Suppose now that $p=0$, so that $\wt{\alpha}=\alpha$, then $\wt{f}(\wt{\alpha})\neq f(\alpha)$ since $\wt{\alpha}\not\in \wt{\cT}$. Therefore $q$ has length 
$\geq 1$ and 
(*) has a factor  $\alpha f(\alpha)g(f(\alpha))$ which is zero in $\La$.

(d)  
So far we have verified that condition (1) of  Definition \ref{def:2.1} holds. 
Condition (3) is also satisfied, from analogous conditions in $\La$. 
We can also see that condition (2) holds: 
For example consider 
$$\wt{\alpha} \wt{f}(\wt{\alpha})\wt{g}(\wt{f}(\wt{\alpha})). \leqno{(**)}$$
If $\wt{\alpha}$ is not in $\wt{\cT}$ then already the product of
the first two terms is zero. Suppose $\wt{\alpha} \in \cT$, then 
(**) is equal to $\alpha f(\alpha)\wt{g}(f(\alpha))$, which has a factor
$\alpha f(\alpha)g(f(\alpha))$ and is zero in $\La$. 
Similarly one obtains the other identity.

\bigskip

(ii)  We investigate when $R$ satisfies the condition (*), that is 
$$\wt{m}_{\wt{\alpha}}\wt{n}_{\wt{\alpha}} \geq 2 \ \mbox{ and } 
\wt{m}_{\wt{\alpha}}\wt{n}_{\wt{\alpha}} 
\geq 3 \mbox{ if $\bar{\wt{\alpha}} \in \wt{\cT}$. } 
$$
Recall $\wt{m}_{\wt{\alpha}} = m_{\alpha}$, hence  if $m_{\alpha}\geq 3$ then this condition holds. 
Assume now that $m_{\alpha}=2$, then the first part of (*) holds.
Suppose that  we have
$m_{\alpha}\wt{n}_{\wt{\alpha}} = 2$,  then we need to show that then $\bar{\wt{\alpha}}$ is not in $\wt{\cT}$.

Write $\wt{\alpha} = \alpha\ldots g^p(\alpha)$, then $B_{\alpha}= \wt{\alpha}^2$, of length $\geq 3$ as an element of $KQ$ (by the assumption on $\La$), and hence
$p\geq 1$.  So we have $t(g^p(\alpha))=i$ but  $s(g^p(\alpha))$ is not in $\Gamma$. 
Assume for a contradiction that
$\bar{\wt{\alpha}}$ is in $\wt{\cT}$, then $\bar{\wt{\alpha}}  = \wt{\bar{\alpha}} = \ba$ and the vertices between $\ba, f(\ba)$ and $f^2(\ba)$ belong to $\Gamma$.
Now, $f^2(\ba) = g^{-1}(\alpha) = g^p(\alpha)$ and therefore $s(g^p(\alpha))$ is in $\Gamma$, a contradiction.
We have shown that when $m_{\alpha}=2$ for an arrow $\alpha$ starting at $i$, the condition (*) holds for $\alpha$.

Assume now that $m_{\alpha}=1$. 
It is possible that  $\wt{n}_{\wt{\alpha}}=1$ so that already the first part of (*) fails. (For example, take $B_{\alpha} = \wt{\alpha}$ of length $\geq 3$ and $s(\alpha)$ is the only vertex
along $B_{\alpha}$ which is in $\Gamma$. This is the exception (*1).)
We continue with $m_{\alpha}=1$, and we assume now 
$\wt{n}_{\wt{\alpha}} = 2$, in this case  the first part of (*) holds.

We  write $B_{\alpha} = (\alpha\ldots g^p(\alpha))(g^{p+1}(\alpha)\ldots g^r(\alpha))$ where $(\alpha \ldots g^p(\alpha)) = \wt{\alpha}$, 
so that we have $\wt{g}(\wt{\alpha}) = (g^{p+1}(\alpha)\ldots g^r(\alpha))$. 
Then 
$i=s(\alpha)$ and $j=s(g^{p+1}(\alpha))$ are the only vertices along the $g$-cycle of $\alpha$
which belong to $\Gamma$.
The condition (*) fails in this case  if and only if $\bar{\wt{\alpha}}$ 
belongs to  $\wt{\cT}$.

We observe that  $\bar{\wt{\alpha}} = \wt{\ba}$, and this belongs to $\wt{\cT}$ if and only 
if all the vertices between  $\ba,  f(\ba)$ and $f^2(\ba)$ belong to
$\Gamma$, that is, each of  $i$ and $t(\ba)$ and $s(f^2(\ba))$ is in $\Gamma$.

We have
$f^2(\ba) = g^{-1}(\alpha) = g^r(\alpha)$, and therefore by the construction
$r=p+1$ and 
the vertex $s(g^r\alpha)$ is what we called $j$. 
In addition we have $t(\ba)$ in $\Gamma$. 
We have arrived at condition (*2).
\end{proof}

\begin{example}\label{ex:2.8} 
	We take the quiver and the weighted surface algebra
	$\La$ as in Example \ref{ex:2.2}, that is we take $\La = H_{\cT}$ with $\cT = Q_1$. The following examples illustrate that the arrows of $\wt{Q}$ need not be a minimal generating set, that is, $\wt{Q}$ may not be the Gabriel quiver of the algebra $e\La e$. 

(a) Let $\Gamma = \{ i\}$. The algebra $R=e\La e$ has the quiver with vertex
	$i$ and two loops, $\wt{\alpha}$ and $\wt{\tau}$. We have $m_{\alpha}=1$ and $\wt{n}_{\wt{\alpha}}=1$ since $\wt{\alpha} = B_{\alpha}$. This is 
	an example for the exception (*1) of Proposition \ref{prop:2.7}.
In fact we also have that $m_{\wt{\tau}}=1$ and $\wt{n}_{\wt{\tau}}=1$. Here 
	$\wt{Q}$ is not the Gabriel quiver of $R$.

(b) Let $\Gamma = \{ i, k, y\}$. Then again $\wt{n}_{\wt{\alpha}}=2$. Now 
	$\wt{\cT} = \{ \tau, \xi, \beta \}$ and $\tau = \ba$.
	The quiver of $R$ is triangular,
	\[
\xymatrix@R=3.pc@C=1.8pc{
%  \xymatrix@C=.8pc{
    k
    \ar@<.35ex>[rr]^{\beta}
    \ar@<.35ex>[rd]^{\omega}
    && i
	\ar@<.35ex>[ll]^{\wt{\alpha}}
    \ar@<.35ex>[ld]^{\tau}
    \\
    & y
    \ar@<.35ex>[lu]^{\xi}
	\ar@<.35ex>[ru]^{\wt{\rho}}
  }
\]
here $\wt{\alpha} = \alpha\eta \delta\gamma$ and $\wt{\rho}  = \rho\sigma$. 
	The permutation $\wt{g}$ is the product of
	three $2$-cycles, 
$$(\wt{\alpha} \ \beta)(\xi \ \omega)(\tau \ \wt{\rho})
	$$
	The arrow $\bar{\wt{\alpha}} = \tau$ is in $\wt{\cT}$ and 
	we have  an example for
	the exception (*2) of Proposition \ref{prop:2.7}.
Note that $m_{\wt{\rho}} = 1$ and $m_{\wt{\alpha}} = 1$. 

(c) Let $\Gamma = \{ i, j, k, y\}$.  The algebra
$R$ has quiver
\[
  \xymatrix@R=1pc{
%   \xymatrix{
   & i \ar[rd]^{\alpha}  \ar@<-.6ex>[ld]^{\wt{g\tau}}  \\
   y    \ar[ru]^{\tau}  \ar@<-.6ex>[rd]^{\omega}
   && j   \ar@(ru,dr)^{\eta}[]  \ar[ld]^{\wt{f\alpha}}   \\
   & k \ar[lu]^{\xi}  \ar[uu]_{\beta}
  }
\]
and $\wt{g} = (\xi \ \omega)(\tau \ \wt{\rho})(\beta \ \alpha \ \eta \ \wt{\delta})$ with multiplicities
$m_{\xi}=2, m_{\tau}=1$ and $m_{\beta}=1$. 
We have
$$\wt{f} = (\omega \ \wt{\rho} \ \alpha \ \wt{\delta} )(\tau \ \xi \ \beta)(\eta)$$
 In this case
the set of distinguished arrow is $\wt{\cT} = \{ \tau, \xi, \beta, \eta\}$.
 We can see directly using identity (2) of Definition \ref{def:2.1} 
that products of arrows in the 4-cycle of $f$ are zero. 

We observe that $m_{\wt{\rho}}\wt{n}_{\wt{\rho}} =2$. and $\bar{\wt{\rho}} = \xi \in \wt{\cT}$, that is 
the multiplcitiy condition is not satisfied. Indeed, we have $s(g^{-1}(\rho))=i\in \Gamma$ and $t(\xi) = k\in \Gamma$ and we have
again an example for the exception (*2) of Proposition \ref{prop:2.7}. 
\end{example}

\bigskip

\section{General hybrid algebras}

We present now the general definition. 
The multiplicity condition (*) in 2.3 is replaced by the weaker requirement
(**). This has the effect that the quiver $Q$ need not be the Gabriel quiver of the algebra, and therefore  we get  many more algebras. 
However now  there are  exceptions for the zero relations, and they are the main reason for much of the work.

We use the notation as in 2.2, 
in particular $\cT$ is a fixed set of triangles (see
2.3). The condition (*) in 2.3 is replaced by the following. 

(**) \ 
We assume $m_{\alpha}n_{\alpha}\geq 2$ for all $\alpha \in Q_1$, except that 
$m_{\alpha}n_{\alpha}=1$ is allowed 
when $\alpha, \ba$ are both not in $\cT$.

Then sometimes an arrow may not be part of the Gabriel quiver, and this motivates our definition of virtual arrows:

\begin{definition} \label{def:3.1} \normalfont Let  $i$ be a vertex, and let
	$\alpha$ be an arrow  starting at $i$. Then $\alpha$ is a {\it virtual} arrow if one of the following holds:\\
	(a) \   $m_{\alpha}n_{\alpha}=1$  and $\alpha, \ba \not\in \cT$; or\\
(b) \ $m_{\alpha}n_{\alpha}=2$ and $\ba \in \cT$.
        That is, $|A_{\alpha}|=1$ and $\ba\in \cT$.
\end{definition}

For the general definition of a hybrid algebra, there are exceptions for the zero relations. To spell these out explicitly, we will use the
term 'critical' as in the following definition.

\begin{definition}\label{def:3.2} \normalfont Let $\alpha$ be an arrow. We say that $\alpha$ is {\it critical}  if
$|A_{\alpha}|=2$ and $\alpha\in \cT$, and moreover $f(\alpha)$ is virtual  (so that $|A_{f(\alpha)}|=1$ and $g(\alpha)\in \cT$).
\end{definition}

In Subsection \ref{ssec:3.1} we 
present diagrams showing the quiver near a virtual arrow, or near a critical arrow.

\begin{definition}\normalfont \label{def:3.3}  Let $(Q, f)$ be a biserial quiver with the data $m_{\bullet}, c_{\bullet}$ as in \ref{ssec:2.2}, and let  $\cT$ be a set of
distinguished  triangles. The  hybrid algebra $H=H_{\cT} = H_{\cT}(Q, f, m_{\bullet}, c_{\bullet} )$, with assumption (**), 
is the algebra $H= KQ/I$ where $I$ is generated by the following elements:

        (1) \
        $\alpha f(\alpha) - c_{\bar{\alpha}}A_{\bar{\alpha}}$\ for  $\alpha \in \cT$  and  $ \alpha f(\alpha)$ for  $\alpha\not\in\cT$.

        (2) \ $\alpha f(\alpha) g(f(\alpha))$
        unless $\alpha, \ba\in \cT$,  and  $\ba$ is either virtual,  or is critical.

        (2') \  $\alpha g(\alpha)f(g(\alpha))$
        unless $\alpha, g(\alpha) \ \in \cT$, and
        $f(\alpha)$ is either virtual,  or is critical.

        (3) \ $c_{\alpha}B_{\alpha} - c_{\bar{\alpha}}B_{\bar{\alpha}}$ for any arrow $\alpha$ of $Q$.

	(4) \ If all arrows of $Q$ are virtual,  then we require $B_{\alpha}\alpha \in I$ and  $\alpha B_{g(\alpha)}\in I$ for each arrow $\alpha$.
\end{definition}

\medskip

If $\cT= Q_1$  and $|Q_0|\geq 2$ this is the same as the definition of a 
weighted surface algebra in \cite{WSA-GV}, but there we did not use the term 'critical'. 
If $\cT = \emptyset$ then the algebra $H_{\cT}$ is special biserial (by (1)), and identities (2) and (2') hold automatically. 
We will  mainly discuss algebras where $\cT\neq \emptyset$.

The details for the  definition of a hybrid algebra are chosen to ensure that
they are precisely the idempotent algebras of weighted surface algebras, up to socle equivalence. 
Furthermore, we require that hybrid algebras are symmetric, and finite-dimensional.
Therefore  a few small algebras need  to be excluded, which actually are the same which were excluded for weighted surface algebras:

\begin{assumption} \label{ass:3.4}\normalfont   
We exclude four algebras, they are not symmetric.\\ 
(1) \  $|Q_0|=2$,  $\cT=Q_1$, with a virtual loop, and
	 the  3-cycle of $g$ has multiplicity $m= 1$ (see \ref{ssec:4.2}(2a)).\\
(2) \ $|Q_0|=3$, $\cT=Q_1$, the singular  algebra with a triangular
	quiver (see \ref{ssec:4.3}(3)), or the singular algebra with a linear quiver (see \ref{ssec:4.4}).\\
	(3) \ $|Q_0| =3$ with a triangular quiver,   $\cT = Q_1$ and $m\equiv 1$ (see \ref{ssec:4.3}(1)).\\
	(4) \ $|Q_0|=6$, $\cT=Q_1$ when $H$ is the singular spherical algebra as in \cite[3.6]{WSA-GV} (see \ref{ssec:4.7}).
\end{assumption}

\bigskip

In \cite[2.7]{WSA-GV},   we had formulated a  slightly different assumption, this is covered by the above 
(modulo minor changes).
One would have liked to have that the Gabriel quiver of $H$ is obtained from $Q$ by removing the
virtual arrows. There is however one exception of a local algebra, which is a hybrid algebra
(it occurs as an idempotent algebra of a weighted surface algebra, see Example \ref{ex:2.8}(a)).

\medskip

\begin{remark}\label{rem:3.5} 
	In the following there will be computations using the
permutations $f$ and $g$, we describe some basic properties. We will use these freely.\\
(1) We always have that  $f^{-1}(\alpha) = g^{-1}(\ba)$. If $\alpha \in \cT$ then $f^{-1}(\alpha) = f^2(\alpha)$ (which may be $\alpha$).\\
(2) \ 
	Assume $i$ is a quaternion  vertex.
	Then 
	we have, exactly as in \cite{WSA, WSA-GV}, 
$$\alpha f(\alpha) f^2(\alpha) = c_{\ba}A_{\ba} f^2(\alpha) = c_{\ba} B_{\ba} = c_{\alpha}B_{\alpha} = \ba f(\ba)f^2(\ba).
	$$
	\end{remark}

\begin{lemma} \label{lem:3.6}Assume $H= KQ/I$ is a hybrid algebra. Then  the Gabriel quiver $Q_H$
of $H$ is obtained from $Q$ by removing the virtual arrows, except when $H$ is local  with two virtual loops.
\end{lemma}

\begin{proof}  Suppose $i$ is a vertex with arrows  
$\alpha, \ba$ starting at $i$. If they are not virtual then they are part of the Gabriel quiver.  As well, if
(say) $\alpha$ is virtual but $\ba$ is not virtual then $\ba$ is part of the Gabriel quiver but $\alpha$ is not.
Suppose now that $\alpha, \ba$ are both virtual.

	(1) Suppose (say) $\alpha$ is a virtual loop and $\ba$ is virtual
	but not a loop.
	Then $\ba$ must be virtual of type (b) as in Definition \ref{def:3.1}, and $m_{\ba}n_{\ba} = 2$ which shows $g(\ba): t(\ba) \to i$, and $\alpha \in \cT$.
	The arrow $f(\alpha)$ starts at $i$, so we have either $f(\alpha) = \alpha$, or $f(\alpha) = \ba$. 
	In the first case we would have $g(\alpha) = \ba= g^2(\ba)$ and $\alpha = g(\ba)$, so that $t(\ba)=i$ and $\ba$ is a loop,  which is not the case. 
	Therefore we can only have
	$f(\alpha) = \ba$, and  since $f^2(\alpha)$ must end at $i$ we have  $f^2(\alpha) = f(\ba): t(\ba) \to i$ and
	it follows that $f(\ba) = g(\ba)$, a contradiction. 
So this cannot happen.

(2) Suppose that $\alpha$  and $\ba$ are virtual but not loops, then they are both in $\cT$ (and they cannot be 
double arrows since  then $g$ would consist of two 2-cycles, and $Q$ would have only two vertices, hence 
the arrows cannot be in 3-cycles of $f$). 
Then $Q$ has a subquiver of the form
\[
  \xymatrix{
%  \xymatrix@C=1pc{
    3
    \ar@<.5ex>[r]^{g(\ba)}
    & i
    \ar@<.5ex>[l]^{\ba}
    \ar@<.5ex>[r]^{\alpha}
    & 2
    \ar@<.5ex>[l]^{g(\alpha)}
  }
%,
\]
with $m_{\alpha}=1=m_{\ba}$. 
By definition of virtual,  $\alpha$ and $\ba$ are in $\cT$,  hence
they must lie in $3$-cycles of $f$. 
Then $f^2(\alpha)$ ends at vertex $1$, so it is either $g(\alpha)$ or $g(\ba)$. 
Since $f(f^2(\alpha)) =\alpha = g(g(\alpha))$ it follows that $f^2(\alpha)\neq g(\alpha)$,  hence it is equal to $g(\ba)$.
Therefore, $f(\alpha)$ must be an arrow $2\to 3$. Similarly $f(\ba)$ is an arrow $3\to 2$.
	That is,  $Q$ is the triangular quiver, with three vertices, and $g$ is a product of $2$-cycles. We have $m_{\alpha} = 1 = m_{\ba}$ and we
	have excluded in Assumption \ref{ass:3.4}(3) that
	$m\equiv 1$. It follows that $m_{f(\alpha)}\geq 2$ and 
	$f(\alpha), f(\ba)$ are not virtual. 
	We will see in Lemma \ref{lem:4.2} that such an  algebra
has finite type,  and  that the Gabriel quiver is obtained by removing the virtual arrows.

(3) Assume both $\alpha, \ba$ are virtual loops. 
First, suppose (say) that  $\alpha$ is in $\cT$, then both $\alpha, \ba$ are
        virtual of type (b). We have $f=(\alpha)(\ba)$ and $g=(\alpha \ \ba)$ with $m_{\alpha}=1$. 
This algebra is dealt with in \ref{ssec:4.1}(2a), and we will see that $H\cong K$.
Hence the Gabriel quiver of $H$ is obtained by removing the virtual arrows.

If $\alpha, \ba$ are not in $\cT$, that is they are virtual of type (a) in 
	Definition \ref{def:3.1}, then 
$m_{\alpha}=m_{\ba}=1$. We see that
$H\cong K[x]/(x^2)$, and that  its Gabriel quiver is not obtained from $Q$ by removing the virtual arrows.
\end{proof}

\begin{corollary}\label{cor:3.7} The only hybrid algebras for which all arrows are virtual are  local algebras \ref{ssec:4.1} (2a) and \ref{ssec:4.1}(1) with $m_{\bullet} \equiv 1$.
\end{corollary}

\begin{proof}
Assume $\alpha$ is virtual of type (a), then $\alpha, \ba$ are not in $\cT$. Since we also assume $\ba$ is virtual 
it must also be of type (a). By (3) of the above proof,  $H$ is as stated.
Suppose now all arrows are virtual of type (b). Then we can  proceed as in part (2) of  the proof of Lemma \ref{lem:3.6},  and
get $H$ is the algebra with triangular quiver and $m\equiv 1$. But this is excluded (see Assumption \ref{ass:3.4}(3)). 
\end{proof}

\bigskip

\subsection{The exceptions in relations (2) and (2')}\label{ssec:3.1}

The exceptions in (2) and (2') of  Definition \ref{def:3.3}  create special cases in various proofs to come.

First we show that there is a unique algebra with two vertices where a critical arrow
occurs in a $g$-cycle with a loop (see Lemma \ref{lem:3.8} below). Otherwise
the exceptions
always arise in specific subquivers of the same kind, for which
we will now fix notation,
to be used later. We write $\zeta_{\alpha} = \alpha f(\alpha) g(f(\alpha))$ and
$\xi_{\alpha} = \alpha g(\alpha) f(g(\alpha))$. We always have $\alpha, \ba \in \cT$, hence all virtual arrows are of type (b). 

We  take care of critical arrows whose $g$-cycle contains a loop.

\begin{lemma}\label{lem:3.8} Assume $\tau$ is critical.

(a) The $g$-cycle of $\tau$ contains a loop if and only if $|Q_0|=2$ and
	$H$ is the algebra in \ref{ssec:4.2}(2c).

 (b) \ Assume the $g$-cycle of $\tau$ does not contain a loop, then $f(\tau)$ cannot be a loop.

\end{lemma}

\begin{proof} Assume $\tau$ is critical, then $|A_{g^i(\tau)}|\neq |A_{f(\tau)}|$ and hence $f(\tau)$ does not
belong to the $g$-cycle of $\tau$. \\
	(a)   For $H$ as in \ref{ssec:4.2}(2c)  one checks 
	directly that the arrow $\tau:= \gamma$ is critical and its $g$-cycle contains a loop.
For the converse, 
assume $\tau$ is critical. If $g(\tau)=\tau$ then $H$ cannot be local (if so then $\tau$ would be in a 2-cycle of $f$). Hence $Q$ contains
 $\xymatrix{
         i \ar@(dl,ul)[]^{\tau} \ar@<+.5ex>[r]^{f(\tau)}
        & j \ar@<+.5ex>[l]^{f^2(\tau)}}
   $
but then since $f(\tau)$ is virtual we have $g(f(\tau)) = f^2(\tau) = f(f(\tau))$ which is a contradiction.
It follows that the $g$-cycle of $\tau$ has length 3 and is a subquiver of $Q$ of the form
$\xymatrix{
        i \ar@(dl,ul)[]^{} \ar@<+.5ex>[r]^{}
        & j \ar@<+.5ex>[l]^{}}. 
   $

Now, $f(\tau)$ is not part of this subquiver but $\tau$ is in $\cT$. It follows that $f(\tau)$ is a loop at $j$ and $\tau$ is the arrow $i\to j$.
	In particular $Q$ has three vertices and  $H$ is the algebra in \ref{ssec:4.2}(2c) with $\gamma$ as the critical arrow.

(b) Suppose $\tau: j\to y$, and assume $f(\tau)$ is a loop. 
Then since $\tau \in \cT$ we must have that $f^2(\tau) : y\to j$. But as well  the arrow $g(\tau)(\neq f(\tau))$ starts at $y$. Since $Q$ is 2-regular, 
we deduce $g(\tau) = f^2(\tau)$ and since $g^3(\tau)=\tau$ it follows that $g^2(\tau)$ is a loop at $j$, a contradiction.
\end{proof}

In the following, we exclude  the algebra  \ref{ssec:4.2}(2c). That is we assume that a critical arrow does not occur in a $g$-cycle with a loop, and that the 
$g$-cycle with a critical arrow has size 3.

\bigskip

\subsubsection{\bf The subquiver around a critical arrow}

{\it We will see that in the exceptional cases}
$$\zeta_{\alpha} \equiv A_{\alpha} \mbox{ and } \xi_{\alpha} \equiv A_{\alpha}.
$$
Let $\tau$ be a critical arrow, in a $g$-cycle of length three, then by definition $\tau$ and $f(\tau)$ belong to 
$\cT$.  In order to  study the paths $\zeta_{\alpha}$ and 
$\xi_{\alpha}$ near $\tau$ in the exceptional cases,  
we also assume that $g^2(\tau)$ belongs to $\cT$.
Then by Lemma \ref{lem:3.8} the quiver near $\tau$ has the following
form
  \[
%  \xymatrix@R=2pc@C=1.5pc{
%  \xymatrix@R=3.5pc@C=1.8pc{
  \xymatrix@R=3.pc@C=1.2pc{
%  \xymatrix@C=.8pc{
    &&& j
    \ar[ld]^{\tau}
    \ar[rrrd]^{}
    \\
    i
    \ar[rrru]^{}
    \ar@<-.5ex>[rr]_(.6){\omega}
    && y
    \ar@<-.5ex>[ll]_(.4){\xi}
    \ar[rd]^{}
    &&&& k
    \ar[llld]^{}
    \\
    &&& x
    \ar[lllu]^{}
    \ar@/_3ex/[uu]^{}
  }
\]
The  permutation $f$ has 3-cycles 
through vertices $j, y, i$ and $y, x, i$ and $j, k, x$. 
At vertex $k$ the quiver there is at least one other arrow, to have
a 2-regular quiver.
We assume that $\tau$ is critical, so that $m_{\tau}=1$ and moreover $\xi =f(\tau)$ is virtual. Since all $f$-cycles shown belong to $\cT$, the arrow $\omega$ is also virtual.

(a) We study the path $\zeta_{\alpha} = \alpha f(\alpha) g(f(\alpha))$ when 
$\ba$ is critical, using the above diagram.
That is we take for $\alpha$ the arrow $j\to k$, so that  $\ba = \tau$. 
Then we have
$$ \zeta_{\alpha} =  c_{\ba}A_{\ba}g(f(\alpha)) 
	=  c_{\tau}\tau g(\tau) f(g(\tau))
	=  c_{\tau}c_{\xi}\tau\xi 
	= c_{\tau}c_{\xi}c_{\alpha} A_{\alpha}.
$$
This must be non-zero since we require that $A_{\alpha}\neq 0$. 
We note that $A_{\alpha} = \alpha \cdot C \cdot f(\alpha)g(f(\alpha))$
where $C$ is a monomial of positive length.

(b) We study the path $\xi_{\alpha} = \alpha g(\alpha) f(g(\alpha))$ when
$f(\alpha)$ is critical, using the above diagram.
Here we take for $\alpha$ the arrow $i\to j$. Then
$$\xi_{\alpha} = \alpha\cdot c_{\tau}A_{\tau}
	= c_{\tau}\alpha \tau g(\tau)
	= c_{\tau}c_{\omega} \omega g(\tau) 
	= c_{\tau}c_{\omega}c_{\alpha}A_{\alpha}
$$
which again must be non-zero.
We note that $A_{\alpha} = \alpha \cdot C \cdot f(g(\alpha))$ where 
again $C$ is a monomial of length $\geq 1$.

\bigskip

{\bf Remark} (a) \ It is not possible that $\alpha$ and $\ba$ are both critical. 
Suppose $\tau = \ba$ and $\alpha: j\to k$ is also critical, then $f(\alpha): k\to x$ is virtual, so there must be an arrow
$x\to k$ and three arrows start at $x$, a contradiction.\\
(b) \ If $\tau$ is critical in a $g$-cycle of length three then in general $g(\tau)$ need not be in $\cT$.

\bigskip

\subsubsection{\bf Subquivers around a virtual arrow}

{\it We will see that in the exceptional cases } 
$$\zeta_{\alpha} \equiv A_{\alpha} \ \mbox{ and } \xi_{\alpha} \equiv A_{\ba}.
$$
(1) \ Assume first that the virtual arrow  is not a loop, then there is 
a pair $\xi, \omega$ of virtual arrows, and 
the quiver contains
\[
%  \xymatrix@R=2pc@C=1.5pc{
%  \xymatrix@R=3.5pc@C=1.8pc{
  \xymatrix@R=3.pc@C=1.8pc{
%  \xymatrix@C=.8pc{
    & i
    \ar[rd]^{}
    \\
    j
    \ar[ru]^{}
    \ar@<-.5ex>[rr]_{\omega}
    && x
    \ar@<-.5ex>[ll]_{\xi}
 \ar[ld]^{}
    \\
   & k
    \ar[lu]
  }
  \]
The arrows shown form two $3$-cycles of $f$, and belong to $\cT$. 
First we assume $|Q_0|>3$, that is $i\neq k$. We assume $\xi, \omega$ are virtual,
then the other arrows in the diagram are not virtual.

\medskip

(a) Consider $\zeta_{\alpha}= \alpha f(\alpha)g(f(\alpha))$ for $\ba$ virtual, then $\ba$ is one of
$\xi$ or $\omega$. 

Consider the case $\ba=\xi$, then we take for $\alpha$ the arrow $x\to k$.
Then
$$\zeta_{\alpha} = c_{\xi}\xi g(f(\alpha))
	=  c_{\xi}c_{\alpha}A_{\alpha} 
$$
and this must be non-zero. 
One can write $A_{\alpha} = \alpha \cdot C$
where $C$ is a monomial of length $\geq 1$. 
When $\ba = \omega$ then we take for $\alpha$ the arrow $x\to i$ and we get
similarly
$$\zeta_{\alpha} = c_{\omega}c_{\alpha}A_{\alpha} 
$$
and we can write $A_{\alpha} = \alpha C$ with $C$ a monomial of length 
$\geq 1$. 

(b) Consider $\xi_{\alpha}= \alpha g(\alpha) f(g(\alpha))$ 
for $f(\alpha)$ virtual, that is
$f(\alpha)= \xi$ or $\omega$.
If $f(\alpha) = \xi$ then we take for $\alpha$ the arrow $i\to j$, and
$$ \xi_{\alpha} = \alpha c_{\xi}\xi
	= c_{\xi}c_{\ba}A_{\ba}$$
and this must be non-zero. We can write $A_{\ba} = C f(g(\alpha))$ where
$C$ is a monomial of positive length.
Suppose $f(\alpha)=\omega$, then we take for $\alpha$ the arrow $k\to x$, and weget
$$\xi_{\alpha} = c_{\omega}c_{\ba} A_{\ba}$$
which must be non-zero. We can write $A_{\ba} = C f(g(\alpha))$ for a monomial
$C$ of positive length.

\bigskip

(2) Now assume $i=k$ so that $|Q_0|=3.$  By  \ref{ssec:4.3}(2) we can assume the multiplicities are not $(m, 1, 1)$
(as this gives   a Nakayama algebra),    and in  \ref{ssec:4.3}(3) we 
deal with multiplicities $(2,2,1)$. This leaves multiplicities $(m_1, m_2, 1)$ where $(m_1, m_2)\neq (2,2)$ and $m_i\geq 2$. 
This case is similar to the above, we omit details.

(3) \ 
Now we consider a virtual loop, and analyze the exceptions.
Here we can use the quiver
\[\xymatrix{
  i \ar@(dl,ul)[]^{\omega} \ar@<+.5ex>[r]^{}
   & k \ar@<+.5ex>[l]^{}}
\]
where $\omega$ is virtual. 
Consider $\zeta_{\alpha} = \alpha f(\alpha)g(f(\alpha))$ when
$\ba$ is virtual using this diagram, that is $\omega = \ba$. We take
for $\alpha$ the arrow $i\to k$.
By assumption $\omega = g(\omega)$ and therefore 
$f$ has cycle $(\omega \ \alpha \ f(\alpha))$. Moreover $g(f(\alpha)) = \alpha$. 
We have
$$\zeta_{\alpha} = c_{\omega}\omega \alpha = c_{\omega}c_{\alpha}A_{\alpha}.
$$
Now 
consider $\xi_{\alpha} = \alpha g(\alpha) f(g(\alpha))$ when $f(\alpha)$ is virtual, using this diagram. That is $f(\alpha) =\omega$. We take for $\alpha$
the arrow $j\to i$. Then $g(\alpha): i\to j$ and $g(\alpha) = f^2(\alpha)$ 
and $f(g(\alpha)) = \alpha$. We have
$$\xi_{\alpha} = \alpha f^2(\alpha)\alpha = \alpha c_{\omega}\omega = 
c_{\omega}c_{\ba}A_{\ba}$$
$\Box$

As in 3.1.1, we can deduce a general description of a path of type $\zeta$ or $\xi$ in a subquiver of the above forms
(allowing also for arrows at $i$ or $k$): 
The following Corollary gives  already Lemma \ref{lem:7.1}

\begin{corollary}\label{cor:3.9}
Consider any path of length three of the form $\zeta_{\sigma}$ or $\xi_{\sigma}$ in the subquiver of 3.1.1 or 3.1.2 shown.\\
(a) If the path does not contain $\xi$ or $\zeta$ then it must be non-zero in $H$.\\
(b) If the path contains $\xi$ or $\zeta$ then it is zero in $H$.
\end{corollary}

Part (a) is implicitly part of the discussion. Part (b) can be seen using the relations (2) and (2') of Definition \ref{def:3.3}.

\bigskip

\subsection{Consistency, bases and dimensions}

This is an update for the case done in Section 2, when virtual arrows are allowed. This may be
found in the Appendix.

\section{Some hybrid algebras with at most three simple modules}

In \cite{WSA} and \cite{WSA-GV} we have excluded small quivers, to avoid 
technical problems obscuring the general structure. However, here
one of the main aims is characterize hybrid algebras as idempotent
algebras of weighted surface algebras. This forces us to  include small algebras
as well.

In this section we consider some  hybrid algebras whose quiver has at most
four vertices. We will mainly discuss algebras where 
 $\cT\neq \emptyset$, and which can have virtual arrows of type (b), for  small
multiplicities.
Note that given $(Q, f)$ and $\cT$, together with $m_{\bullet}, c_{\bullet}$, the algebra is completely determined, and we will usually not
write down relations explicitly.

\subsection{Local algebras}\label{ssec:4.1}
Here $Q$ consists of one vertex and two loops, denoted by $\alpha$ and $\beta$.
There are two possibilities for $f$ and $g$, and if $f$ is the identity 
permutation there are three possibilities, depending on $\cT$.

\medskip

(1) Consider an algebra where 
$f\ = \ (\alpha \ \beta)$ and  $g = (\alpha)(\beta)$, then we must have $\cT=\emptyset$. 
We may assume $m_{\alpha}\geq m_{\beta}$. \\
{\it If $m_{\beta}=1$ then $H\cong K[x]/(x^{m_{\alpha}+1})$. Otherwise
it is an algebra of dihedral type as in \cite[III.1(a)]{E1}:} 

\medskip

The relations are: 
$$\alpha\beta = 0 = \beta\alpha, \ \ c_{\alpha}B_{\alpha} = c_{\alpha}\alpha^{m_{\alpha}} = c_{\beta}\beta^{m_{\beta}} = c_{\beta}B_{\beta}.
$$
If $m_{\beta}=1$ so that $\beta$ is virtual (of type (a) of Definition \ref{def:3.1}),  then 
$H\cong K[x]/(x^{m_{\alpha}+1})$.
This also holds when $m_{\alpha}=1$; in this case   the Gabriel quiver
of $H$ is not obtained from $Q$ by removing the virtual arrows (see also Lemma \ref{lem:3.6}). If $m_{\beta}\geq 2$  then $H$ is special biserial, of infinite type and is a (commutative) algebra of dihedral type, as defined in
\cite[III.1(a)]{E1}.

\medskip

(2)  Consider hybrid algebras 
where $f=(\alpha)(\beta)$ and $g=(\alpha \ \beta)$, so $m_{\alpha}n_{\alpha}\geq 2$.

(2a) Assume first that  $\cT= Q_1$.
{\it If $m_{\alpha}=1$ then $H\cong K$, and if  
 $m_{\alpha} \geq 2$ then $H$ is an algebra as in \cite[III.1(e)]{E1} of quaternion type:} 

\medskip

Assume  that $m_{\alpha} =1$, we may assume that
$c_{\alpha}=1$.
The relations are
$$\alpha^2 = A_{\beta} = \beta \ \mbox{ and}  \ \beta^2 = A_{\alpha} = \alpha,$$
that is, both arrows are virtual.
By condition (4) of Definition \ref{def:3.3} we have  that $B_{\alpha}\alpha=0 = \alpha B_{g(\alpha)}$.
Relation (3) gives  $B_{\alpha} = \alpha\beta = B_{\beta} = \beta\alpha$. 
and hence $\alpha^2\beta = 0$ and therefore 
$$0 = \alpha^4 = \beta^2 = \alpha
$$
and similarly $\beta=0$.  We have shown that  $H\cong K$. On the other hand, when 
$m_{\alpha}\geq 2$ then we see directly that we get an algebra of quaternion type, as in \cite[III.1.(e)]{E1}.
The algebras where $H\cong K$ cannot occur as an idempotent algebra of a WSA $\La$, since $e_i\La e_i$
has at least two independent elements:  the idempotent $e_i$ and the generator of the socle of $e_i\La$.

(2b) Assume $\cT=\{ \beta \}$.  
{\it If $m_{\alpha}=1$ then $H\cong K[x]/(x^4)$. Otherwise $H$ is an algebra
as in \cite[III.1(d)]{E1} (of semidihedral type):}

We may assume $c_{\alpha}=1$, and we 
have the relations
$$\beta^2 = A_{\alpha} = (\alpha\beta)^{m_{\alpha}-1}\alpha, \ \ \alpha^2 = 0.
$$
If $m_{\alpha}  \geq 2$, this gives precisely the algebras in  \cite[III.1(d)]{E1}. 
Suppose 
$m_{\alpha}= 1$ so that the arrow  $\alpha$ is virtual. 
Then we see $\beta^3 \equiv B_{\alpha}$ and $\beta^4=0$ which shows
that $H$ is isomorphic to $K[x]/(x^4)$.  In this case the Gabriel quiver is obtained from $Q$ by removing the virtual arrows.

\bigskip

(2c) Assume $\cT=\emptyset$. {\it Then $H$ is an algebra as in \cite[III.1(b)]{E1}.  For $m_{\alpha}=1$ it is four-dimensional commutative:}
This is seen directly from the relations.
$\Box$

\bigskip

\subsection{Hybrid algebras with two simple modules}\label{ssec:4.2}

Let $H$ be a hybrid algebra with two simple modules, then $H=KQ/I$ where the  quiver
$Q$  is of the form
\[
  \xymatrix{
%  \xymatrix@C=1pc{
    1
    \ar@(ld,ul)^{\alpha}[]
    \ar@<.5ex>[r]^{\beta}
    & 2
    \ar@<.5ex>[l]^{\gamma}
    \ar@(ru,dr)^{\sigma}[].
  }
%,
\]
We consider only the 
cycle structures of $f, g$ for which  $\cT$ can be non-empty and the algebra can have virtual arrows of type (b).

(1) Consider algebras with 
$$
f=(\alpha)(\beta \gamma)(\sigma) \ \mbox{with  } \ g = (\alpha \ \beta \ \sigma \ \gamma).$$
Suppose $\cT\neq \emptyset$, then 
$\cT$ consists of one or two loops, and there are no virtual arrows. 
The algebras look similar to algebras of semidihedral type in \cite{E1}, however they have always singular Cartan matrices, which was
excluded for semidihedral type. 

\bigskip

(2) Consider algebras where
$$
f = (\alpha \ \beta \ \gamma)(\sigma) \ \ \mbox{ with } \ \
g = (\alpha) (\beta \ \sigma \ \gamma).
$$
For hybrid algebras with $\cT\neq \emptyset$, the possibilities for
 $\cT$ are  either $Q_1$, \ or  $\cT= \{ \sigma\}$, \ or
        $\cT = \{ \alpha, \beta, \gamma\}$.

\bigskip

4.2(2a) \ {\bf The case $\cT = Q_1$ and $(t, m)=(2, 1)$ }. 
This is excluded in Assumption \ref{ass:3.4}(1).
In \cite{WSA-GV} it was  excluded because the algebras  appeared to be of finite type.
However the argument was based on the incorrect relations.
Here we review this algebra, with amended relations.

We may take $c_{\bullet} = (1, c)$. Note that $\alpha$ is virtual and $\gamma$ is critical.
The associated hybrid algebra
is given by the relations
\begin{align*}
\beta\gamma &= \alpha, &
 \gamma \alpha &= c \sigma \gamma,
 &
 \alpha\beta &=  c \beta \sigma,
 &
        \sigma^2 &=  c\gamma\beta, \\
 \alpha\beta\sigma &=  0,
 &
 \gamma \alpha^2 &=  0,
 &
\sigma\gamma\alpha&= 0,  & \alpha^2\beta &=0
 \end{align*}

These imply that the algebra is not symmetric. Alternatively, there
is a quick way to get a contradiction. Namely
$$0=\beta\sigma^2 \equiv \beta\gamma\beta = \alpha\beta$$
and $\alpha \not\in \cT$.

\bigskip

 4.2 (2b) {\bf The case  $\cT = Q_1$ and
        $t=3$ and   $m=1$.}  \
This was  dealt with in \cite[Example 3.1(1)]{WSA-GV}, the algebra is called {\it disc algebra}, and is
denoted by $D(\lambda)$. Viewed in the context of  periodicity,
it has a
singular version: when the scalar parameter $\lambda =1$ it is not periodic,
but it is a hybrid algebra.
In that case ${\rm rad}(e_1H)/S_1\cong {\rm rad}(e_2H)/S_2$ and is indecomposable,  and the simple modules belong to an Auslander-Reiten component
of type $D$. The algebra is
 of semidihedral type, part of the family $SD(2\cB)_3$ in \cite{E1}
and it is a hybrid algebra.

\bigskip

4.2 (2c) \ {\bf Algebras with $\cT= \{\alpha, \beta, \gamma\}$ and $(t, m)  = (2, 1)$.} This is the only algebra where the
$g$-cycle of a critical arrow has a loop (see Lemma \ref{lem:3.8}). However the algebra is seen below to be special biserial and we do not
have to consider it further.
The arrow  $\alpha$ is virtual, and $\gamma$ is critical.
We may take $c_{\alpha}=1$, and we set $c_{\beta}=c$. Then the associated hybrid algebra is given by the relations:
\begin{align*}\beta\gamma&= \alpha, & \gamma\alpha&= c\sigma\gamma= \gamma\beta\gamma, & \alpha\beta&= c\beta\sigma=\beta\gamma\beta, &
        \sigma^2&= 0, \cr
        \beta\gamma\beta\sigma&=0, & (\gamma\beta)^2\gamma&=0, & (\beta\gamma)^2\beta&=0, & \sigma^2\gamma&=0, & \beta\sigma^2&=0, & \sigma\gamma\beta\gamma&=0.
\end{align*}
Note that  $\gamma\beta\sigma = \sigma\gamma\beta = c^{-1}(\gamma\beta)^2$
and $B_{\beta}J=0=B_{\gamma}J$.

\begin{lemma} \label{lem:4.1} The algebra $H$ is special biserial. More precise, let
        $\sigma':= (c\sigma - \gamma\beta)$. Then $\sigma'\gamma=0$ and $\beta\sigma'=0$. \\
Then  $H$ has presentation
$k\wt{Q}/\wt{I}$ where
$\wt{Q}$ is the quiver
$$
 \xymatrix{
  1  \ar@<+.5ex>[r]^{\beta}
  & 2 \ar@<+.5ex>[l]^{\gamma} \ar@(ur,dr)[]^{\sigma}
 }
      $$
and $\wt{I}= \langle \sigma\gamma, \ \beta\sigma,  \   \sigma^2 -  (\gamma\beta)^2\rangle $. 
\end{lemma}

\bigskip

{\it Proof } \ Rewriting the relations gives that $\sigma'\gamma=0$ and $\beta\sigma'=0$. Note that $\sigma'$ may be taken
as an arrow. 
We have  $\sigma'\sigma =  c\sigma^2 - \gamma\beta\sigma$ and it is non-zero in the socle of $e_2A$.
We have
$\sigma'\gamma\beta = c(\sigma\gamma\beta - (\gamma\beta)^2) = 0$, hence
$$(\sigma' )^2 =   -(\gamma\beta)^2$$
We may rescale $\sigma'$ and then obtain the presentation as stated.
$\Box$

\medskip

\parindent=0pt One may introduce a virtual loop at $1$, which gives 
  a presentation of a hybrid algebra.
             \bigskip

4.2 (2d) \ {\bf Algebras with  $\cT = \{ \sigma\}$ and
$t=m=1$. }

Here $\alpha$ is virtual of type (a) (note that $\alpha$ and $\ba = \beta$ are not in $\cT$).
We can take $c_{\alpha}=1$ and we set $c_{\beta}=c$.
Then the relations are
\begin{align*}
        \alpha\beta&=0, & \beta\gamma&=0, & \gamma\alpha&=0, &
        \sigma^2\gamma&=0,  \cr  \beta\sigma^2&=0, &
        \sigma^2&=c\gamma\beta, & \alpha &= c(\beta\sigma\gamma), & cB_{\gamma} &= cB_{\sigma}.
\end{align*}
This
algebra  occurs in (3.6) of \cite{Sk}, with a slightly different presentation. It is an algebra of finite (Dynkin) type
 $\mathbb{D}$.

\bigskip

We consider now some algebras with
three simple modules.

In total there are five possible quivers for which $f$ has at least one 3-cycle.  We will discuss algebras with
'triangular' and 'linear' quiver in some detail first, and will briefly consider the other three later.

\medskip

\subsection{Algebras with  triangular quiver}\label{ssec:4.3}
Let $Q$ be the quiver
\label{ex:3.3}
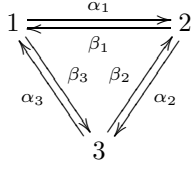
\[
 \begin{tabular}{c}
%  \xymatrix@R=2pc@C=1.5pc{
%  \xymatrix@R=3.5pc@C=1.8pc{
%  \xymatrix@R=3.pc@C=1.8pc{
%  \xymatrix@C=.8pc{
 %   1
  %  \ar@<.35ex>[rr]^{\alpha_1}
  %  \ar@<-.35ex>[rr]_{\beta_1}
%  && 2
  %  \ar@<.35ex>[ld]^{\alpha_2}
       %  \ar@<-.35ex>[ld]_{\beta_2}
 %  \\
  %  & 3
                  %  \ar@<.35ex>[lu]^{\alpha_3}
  %  \ar@<-.35ex>[lu]_{\beta_3}
 % }
%  \xymatrix@R=2pc@C=1.5pc{
%  \xymatrix@R=3.5pc@C=1.8pc{
           \xymatrix@R=3.pc@C=1.8pc{
%  \xymatrix@C=.8pc{
    1
    \ar@<.35ex>[rr]^{\alpha_1}
    \ar@<.35ex>[rd]^{\beta_3}
    && 2
    \ar@<.35ex>[ll]^{\beta_1}
    \ar@<.35ex>[ld]^{\alpha_2}
    \\
    & 3
    \ar@<.35ex>[lu]^{\alpha_3}
   \ar@<.35ex>[ru]^{\beta_2}
  }
  \end{tabular}
   \]
The only cycle structure for which   $\cT$ can be non-empty is given by
  $f \ =  (\alpha_1\ \alpha_2\ \alpha_3)(\beta_1\ \beta_3\ \beta_2)$, so that
    $g= (\alpha_1\ \beta_1)(\alpha_2 \ \beta_2)(\alpha_3 \ \beta_3)$. We write $m_i = m_{\alpha_1}$ and $c_i = c_{\alpha_i}$. 
      
   \medskip

4.3 (1)  {\bf   Algebras with   $\cT=Q_1$ and $m_{\bullet} = (1,1,1)$. }
 Such an algebra is   excluded in Assumption
\ref{ass:3.4}(3). It was  excluded in \cite[4.4]{WSA-GV}, though the argument was not  correct.  
The algebra is given by the  relations
$$\alpha_i\alpha_{i+1} = c_{i+2}\beta_{i+2}, \ \ \beta_i\beta_{i-1} = c_{i-2}\alpha_{i-2}
$$
(indices modulo 3). As well $B_{\alpha_i}= \alpha_i\beta_i \equiv B_{\beta_{i-1}} = \beta_{i-1}\alpha_{i-1}$, and there are no
zero relations of types (2) or (2').  
 We observe that
$$\alpha_1\beta_1 \equiv \beta_3\beta_2\alpha_2\alpha_3 \equiv \beta_3\alpha_3\beta_3\alpha_3 \equiv (\alpha_1\beta_1)^2 = 0$$
and this is zero by condition (4) of Definition \ref{def:3.3}.
Similarly all paths $\alpha_i\beta_i$ and $\beta_i\alpha_i$ are
zero, and then any cyclic path of positive length is zero in the algebra.
Therefore the algebra is not symmetric.

\medskip

 4.3 (2) \ {\bf Algebras with  $\cT=Q_1$ and $m_{\bullet} = (m, 1, 1)$ and $m\geq 2$.}\\
Such an 
algebra was   excluded in \cite[4.4]{WSA-GV},
as it was said to be  not finite-dimensional.
However this is not correct, it has  even finite type, as we will now show.
Note also that the Gabriel quiver is obtained by removing the virtual arrows.

\begin{lemma}\label{lem:4.2}  With these conditions, $H$ has finite type, it is isomorphic to the direct sum of a Nakayama algebra
$$KQ/\langle (\alpha\beta)^{m-1}\alpha, (\beta\alpha)^{m-1}\beta\rangle
$$
with a copy of $K$, where $Q$ is the quiver
$ \xymatrix{
  1  \ar@<+.5ex>[r]^{\alpha}
   & 2 \ar@<+.5ex>[l]^{\beta}
 }.
$
\end{lemma}

{\it Proof }
The relations are as follows.
\begin{align*} \alpha_1\alpha_2=& c_3\beta_3, & \alpha_2\alpha_3 =& c_1A_{\beta_1}, & \ \alpha_3\alpha_1 =& c_2\beta_2\cr
\beta_1\beta_3=& c_2\alpha_2, & \beta_3\beta_2=& c_1A_{\alpha_1}, & \beta_2\beta_1 =& c_3\alpha_3, \cr
\alpha_2\alpha_3\beta_3 &=0, & \beta_3\beta_2\alpha_2&=0, &
\alpha_3\beta_3\beta_2&=0, & \beta_2\alpha_2\alpha_3&=0.
\end{align*}
Moreover we have the consequences
\begin{align*}c_1B_{\alpha_1} = c_3B_{\beta_3}, &&  c_1B_{\beta_1} = c_2B_{\alpha_2}, && c_2B_{\beta_2} = c_3B_{\alpha_3}
\end{align*}

(1) Starting with the relation $0=\alpha_2\alpha_3\beta_3 (= \alpha_2B_{\alpha_3})$ we show that $\beta_1\alpha_1\alpha_2=0$: Namely
$$0=\alpha_2B_{\alpha_3} \equiv \alpha_2B_{\beta_2} = B_{\alpha_2}\alpha_2 \equiv B_{\beta_1}\alpha_2 = (\beta_1\alpha_1)^m\alpha_2.
$$
Next we have
$$(\beta_1\alpha_1)^m\alpha_2 = (\beta_1\alpha_1)^{m-1}\beta_1\alpha_1\alpha_2 \equiv (\beta_1\alpha_1)^{m-1}\beta_1\beta_3 \equiv (\beta_1\alpha_1)^{m-1}\alpha_2.
$$
Repeating this reduction gives $\beta_1\alpha_1\alpha_2=0$ and then $\alpha_2=0$. 
Similarly we have  $0=\beta_3=\beta_2=\alpha_3$.  Hence the algebra
has a direct summand spanned by $e_3$ which is  isomorphic to $K$. Furthermore, from the relations we have
$A_{\beta_1}=0$ and  $A_{\alpha_1}=0$, and there are no further restrictions.
This shows  that the subalgebra generated by $e_1, e_2$ and $\alpha_1, \beta_1$ is the Nakayama algebra as stated.
$\Box$

\bigskip
4.3 (3) \  {\bf Algebras with  $\cT=Q_1$ and  $m_{\bullet} = (2, 2, 1)$.}

 They  are  called   {\it triangle algebras, } as discussed in  \cite[Example 3.3 (1)]{WSA-GV},  and  denoted by $T(\lambda)$,  where $c_{\bullet} = (\lambda, 1, 1)$.
 The algebra with $\lambda=1$  is not symmetric, as it was shown in \cite[3.3]{WSA-GV}, and therefore it is excluded in Assumption \ref{ass:3.4}(2).

\bigskip

4.3 (4) \ {\bf Algebras with  $\cT= \{ \alpha_1, \alpha_2, \alpha_3\}$ and $m_{\bullet} = (1, 1, 1)$. }\ 
In this case, the arrows $\beta_i$ are virtual, 
and  the algebra $H_{\cT}$
is a Nakayama algebra of finite type: The relations are
$$\alpha_i\alpha_{i+1} = c_{i+2}\beta_{i+2} \ \mbox{and} \  \beta_i\beta_{i-1}=0. 
$$
There are no exceptions to the zero relations in (2) and (2') since for any arrow $\alpha$ we have  $\alpha\not\in \cT$ or $\ba\not\in \cT$, and 
$\alpha\not\in \cT$ or $g(\alpha)\not\in \cT$. It is straightforward to check 
that  $H$ 
is the Nakayama algebra where the quiver is cyclic with three
vertices, and where all paths of length four are zero
in the algebra.
%This occurs as an idempotent algebra, see Example \ref{ex:5.3}(e).

\bigskip

\subsection{Algebras with linear quiver}\label{ssec:4.4}

Consider algebras whose quiver is of the form
\[
  \xymatrix{
%  \xymatrix@C=1pc{
    1
%    \ar `ld_u[] `_rd[]^{\alpha} []
    \ar@(ld,ul)^{\alpha}[]
    \ar@<.5ex>[r]^{\beta}
    & 2
    \ar@<.5ex>[l]^{\gamma}
%    \ar `ru_d[] `_lu[]^{\eta} [] &
    \ar@<.5ex>[r]^{\sigma}
    & 3
    \ar@<.5ex>[l]^{\delta}
%    \ar `ru_d[] `_lu[]^{\xi} [] &
    \ar@(ru,dr)^{\eta}[]
  }
%,
\]
  To have that $\cT\neq \emptyset$ containing some
virtual arrows of type (b), we have
two possibilities for the permutations $f$ and $g$:
$$ f = \ (\alpha \ \beta \ \gamma)(\sigma \ \eta \ \delta) \ \mbox{ and } \ \  g= \ (\alpha)(\beta \ \sigma \ \delta \ \gamma)(\eta), \ \ \mbox{ or } \ \ 
f= \ (\alpha \ \beta\ \gamma)(\sigma \ \delta)(\eta) \ \  \mbox{ and } \ \  g= \ (\alpha)(\beta \ \sigma \ \eta \ \delta \ \gamma).
$$
For most of the hybrid algebras with these cycle structures, virtual arrows do not lead to special cases.
We only discuss algebras with the first cycle structure  and where $m_{\bullet} = (2, 1, 2)$.
This has been considered in Example 3.4 of \cite{WSA-GV}. It is shown that we may assume
$c_{\bullet}= (1, \lambda, 1)$ ,  the algebra
is called $\Sigma(\lambda)$. Furthermore, it is proved (in Lemma 3.5 of \cite{WSA-GV}) that
$\Sigma(\lambda)$ is isomorphic to the triangular algebra $T(\lambda^{-2})$ introduced in \ref{ssec:4.3}. In particular this implies that
we must exclude $\lambda = \pm 1$, since then the algebra is not symmetric.
We refer to this as a singular algebras, which are excluded in Assumption \ref{ass:3.4}(2).

\bigskip

\subsection{Three other quivers with three vertices}\label{ssec:4.5}

The following three quivers also have each at least one 3-cycle of $f$ which may or may not belong to $\cT$:

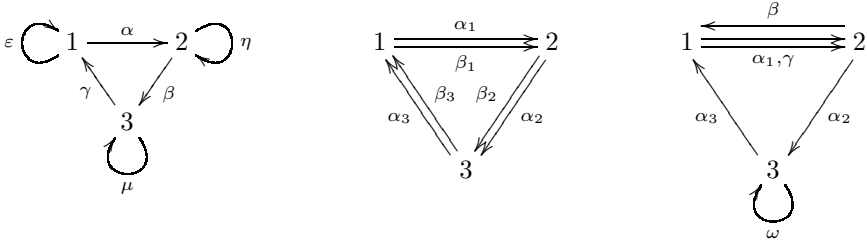
\[
%  \xymatrix@C=.8pc
  \xymatrix@C=.8pc@R=1.5pc
  {
     1 \ar@(dl,ul)[]^{\varepsilon} \ar[rr]^{\alpha} &&
     2 \ar@(ur,dr)[]^{\eta} \ar[ld]^{\beta} \\
     & 3 \ar@(dr,dl)[]^{\mu} \ar[lu]^{\gamma}}
\qquad
\qquad
%  \xymatrix@R=2pc@C=1.5pc{
%  \xymatrix@R=3.5pc@C=1.8pc{
  \xymatrix@R=3.pc@C=1.8pc{
%  \xymatrix@C=.8pc{
    1
    \ar@<.35ex>[rr]^{\alpha_1}
    \ar@<-.35ex>[rr]_{\beta_1}
  && 2
    \ar@<.35ex>[ld]^{\alpha_2}
    \ar@<-.35ex>[ld]_{\beta_2}
    \\
    & 3
    \ar@<.35ex>[lu]^{\alpha_3}
    \ar@<-.35ex>[lu]_{\beta_3}
  }
\qquad
\qquad
%\]
%\[
%  \xymatrix@R=2pc@C=1.5pc{
%  \xymatrix@R=3.5pc@C=1.8pc{
     \xymatrix@R=3.pc@C=1.8pc{
%  \xymatrix@C=.8pc{
    1
    \ar@<.35ex>[rr]^{}
    \ar@<-.35ex>[rr]_{\alpha_1, \gamma}
        && 2
    \ar@<.35ex>[ld]^{\alpha_2}
     \ar@<-1.5ex>[ll]_{\beta}
    \\
    & 3
    \ar@<.35ex>[lu]^{\alpha_3}
    \ar@(dr, dl)[]^{\omega}
      }
\]
For  the first two quivers, there are no virtual arrows of type (b) since there is just one $g$-orbit of size $6$.
Consider the third quiver when $f = (\alpha_1 \ \alpha_2 \ \alpha_3)(\omega)(\beta \ \gamma)$.  Then
$g = (\alpha_1 \ \beta)(\gamma \ \alpha_2 \ \omega \ \alpha_3)$. We consider the case when $m_{\bullet}\equiv 1$, then
the arrow $\beta$ is virtual if $\cT$ contains $\{ \alpha_1, \alpha_2, \alpha_3\}$. However this does not create complications:
If $\beta$ is virtual then relations $\gamma\alpha_2\alpha_3$ and $\alpha_2\alpha_3\gamma$ are excluded in (2), (2') of Definition
\ref{def:3.3}. In this case they are, up to non-zero scalars, equal to $\gamma\beta$ and $\beta\gamma$,  which are zero since $f$ has the cycle $(\beta \ \gamma)$.
We note that the algebras are  not of semidihedral type, as the Cartan matrices 
are singular.

\subsection{An exceptional algebra with four simple modules}\label{ssec:4.6}

Let $H=KQ/I$ where $Q$ is the quiver
\[
%  \xymatrix@R=2pc@C=1.5pc{
%  \xymatrix@R=3.5pc@C=1.8pc{
  \xymatrix@R=3.pc@C=1.2pc{
%  \xymatrix@C=.8pc{
    &&& 1
        \ar[ld]^{\alpha}
    \ar@/^5ex/[dd]^{\ba}
    \\
    4
    \ar[rrru]^{\delta}
    \ar@<-.5ex>[rr]_(.6){\eta}
    && 2
    \ar@<-.5ex>[ll]_(.4){\xi}
    \ar[rd]^{\gamma}
    &&&& 
     \\
    &&& 3
    \ar[lllu]^{\sigma}
    \ar@/_3ex/[uu]^{\beta}
  }
\]
with  $f = \ (\ba \ \beta)(\alpha\ \ \xi\ \delta)(\gamma\ \sigma\ \eta)$, and hence 
$g =  (\alpha\  \gamma \ \beta)(\ba \ \sigma \ \delta)(\xi\ \eta)$.
Moreover, we take $m_{\bullet} = 1$  and $c_{\alpha}=c$ and $c_{\beta}=c_{\xi}=1$.
Let   $\mathcal{T} = \{ \alpha, \ \xi,  \ \delta,  \ \gamma \ \sigma \ \eta\}$, and let $H$ be the hybrid algebra defined by  these data.  
Then 
$\xi$ and $\eta$ are virtual arrows and 
the Gabriel quiver $Q_H$ is obtained by removing $\xi, \eta$.

\begin{lemma}\label{lem:4.3}
The algebra $H$ is special biserial. Let $\bar{Q}$ be the quiver obtained from $Q$ by removing $\xi$ and $\eta$, and adding virtual loops 
$\ve, \rho$ of type (a). Then $H$ has a hybrid algebra presentation with this quiver, and  with $\bar{\cT} =\emptyset$, defined by   by the data 
		data
	$$\bar{f} = (\delta \ \ba' \ \sigma \ \ve)(\alpha \ \rho \ \gamma \ \beta'), \ \ \bar{g} = (\delta \ \alpha \ \gamma \  \sigma)(\ba' \ \beta')(\rho)(\ve)
$$
	with multiplicity $\equiv 1$ and parameter function 
	$\equiv 1$. The loops $\ve, \rho$ are
	virtual of type (a).
\end{lemma}

\begin{proof} Starting with the given presentation, we replace $\beta$ by
	$\beta':= \sigma\delta - c\beta$, then $\beta'\alpha=0$ and
	$\gamma\beta'=0$. We also replace $\ba$ by $\ba':= \alpha\gamma - \ba$, and then $\ba'\sigma = 0$ and $\delta\ba'=0$. We take $\ve$ to be the socle
	monomial $\delta\alpha\gamma\sigma$, and we take $\rho$ to be
	the socle monomial $\gamma\sigma\delta\alpha$.
	Then it is straightforward to show that the algebra has
	the stated presentation.
\end{proof}

\bigskip

\subsection{Singular algebras}\label{ssec:4.7}
In addition to the singular disk, and triangle algebra as we have discussed above, there are two further algebras which
were called singular in \cite{WSA} and \cite{WSA-GV}.
Recall from \cite{WSA} Example 6.1 the {\it tetrahedral algebras}. This family contains one algebra, with certain parameters,
which is not periodic, and therefore it was called singular in that context. However, it is a hybrid algebra.

Furthermore, in Example 3.6 of \cite{WSA-GV} we have discussed
{\it spherical algebras}, denoted by $S(\lambda)$ for $\lambda \in K^*$. The quiver has six vertices, and with the smallest multiplicities the algebra has
four virtual arrows. When $\lambda = 1$, it is not symmetric and is therefore excluded in Assumption \ref{ssec:4.3}(4).

\bigskip

\section{Hybrid algebras as idempotent algebras of weighted
surface algebras}

In the first part of this section we will prove that
for a weighted surface algebra $\La$ and  an idempotent $e$ of $\La$,  every block of $e\La e$ is a hybrid algebra.
In the second part of this section  we will show that
every hybrid algebra with $\cT \neq Q_1$ occurs in this way.
The second part generalizes the main results of \cite{BGA},  which dealt with the hybrid algebras where $\cT = \emptyset$, that is, the Brauer graph algebras.
Note that we start with a weighted surface algebra, which is not a socle deformation.
\bigskip

\begin{theorem}\label{thm:5.1}
Assume $\La$ is a weighted surface algebra and let  $e\in \La$ be an idempotent. Then
each block of the algebra
$e\La e$ is a hybrid algebra.
\end{theorem}

\bigskip

\begin{proof} 
We fix a weighted surface algebra $\La$, and we proceed as in the proof of Proposition \ref{prop:2.7}. 
  By general theory, we may assume that
$e = \sum_{i\in \Gamma} e_i$ with $\Gamma$ a subset of the vertices
of $Q$, and we set $R=e\La e$, and  we may assume that $e$ is not the identity of $\La$. 
We take the  quiver $\wt{Q}$ with vertices labelled by $\Gamma$. For $\alpha \in Q_1$,
	let $\wt{\alpha}$ be the shortest path in $Q$ along the $g$-cycle
	of $\alpha$ starting with $\alpha$ and ending at
	some vertex in $\Gamma$.  We take the set $\wt{Q}_1$ of these $\wt{\alpha}$ as arrows for $\wt{Q}$, it is 
	a generating set for $R$, and we have
	a surjective algebra map $\psi: K\wt{Q} \to R$.
	As in \ref{prop:2.7}, the quiver $\wt{Q}$ is 2-regular. 
	When $\wt{\alpha}=\alpha$ then we write for simplicity $\alpha$.
We define the permutation $\wt{f}$, and the distinguished
	set $\wt{\cT}$ of triangles, as in Proposition \ref{prop:2.7}.
The cycles of the associated permutation $\wt{g}$ are obtained
	from the cycles of $g$ by replacing
	$\alpha, g(\alpha), \ldots, g^{p}(\alpha)$ by $\wt{\alpha}$. 
We take the multiplicity and  parameter functions 
	as for $\La$.  
	Then we may write down elements $B_{\wt{\alpha}}$ of $R$
	for each arrow $\alpha$, and it is clear that these satisfy
	identity (3) of Definition \ref{def:3.3}.
	As well we have
	elements $A_{\wt{\alpha}}$ such that $A_{\wt{\alpha}}\wt{\gamma} = B_{\wt{\alpha}}$ where $\wt{\gamma}$ is the last arrow in $B_{\wt{\alpha}}$. 
 Furthermore, the exceptions in relations (2) and (2') 
	occur precisely when the arrows $\alpha, \ba$ (or $\alpha, g(\alpha)$) are in $\wt{\cT}$.

\medskip
We will show  that the arrows in $\wt{Q}_1$ satisfy the
	identity (1) of Definition \ref{def:3.3}. For
	the arrows in $\wt{\cT}$, this follows directly from identity (1) for
	$\La$. 
	Let $\wt{\alpha}$ be an arrow of $\wt{Q}$ which is not in $\wt{\cT}$, 
	and let  $p:= \wt{\alpha}\wt{f}(\wt{\alpha})$ We must show that this
	is zero in $R$,  
	(possibly after some adjusting), or possibly that  it is a scalar multiple of a socle element, ie we have a socle deformation.

\medskip
	Since $\wt{\alpha}$ is not in $\wt{\cT}$, we know that
	$p$ has length $|p| \geq 3$ as a path in $Q$.
	If $|p| \geq 5$ then it is zero in $\La$, 
	this follows from Lemma \ref{lem:7.5}. 
Suppose now that $p$ is non-zero, then we must have
$|p|=3$ or $|p|=4$. 
	For the following we exclude the algebras \ref{ssec:4.3}(2) (this can be done by hand, using Lemma \ref{lem:4.2}. Furthermore we exclude
	\ref{ssec:4.2} and \ref{ssec:4.3}(3), they will be considered
	below in \ref{ssec:5.1}.

\bigskip

(a) Assume first that $|p|=3$, then $p$  is of the form
	$\zeta_{\alpha}$ of $\xi_{\alpha}$, near a critical or virtual arrow.
	We start with $p$ near a critical arrow.\\
(a1)  {\it Assume  $p = \zeta_{\alpha} = \alpha f(\alpha)g(f(\alpha))$ and $\ba$ is critical.} 
That is we have
	$\wt{\alpha} = \alpha$ and
	$\wt{f}(\wt{\alpha})= f(\alpha)g(f(\alpha))$.We use diagram 3.1.1, and set $\tau = \ba$ so that
	$\alpha: j\to k$. In this case $\Gamma$ contains vertices $j, k, i$ but $\Gamma$ does not contain $x$.	Let $\beta = f(\alpha)$. 
The cycle of $\wt{f}$ containing $\wt{\alpha}$ is
$$(\alpha \ \wt{\beta} \ \wt{\omega}  \gamma \ \wt{\ba})
$$ where $\gamma: i\to j$. 
Note that $\wt{\omega} = B_{\omega}$ and $\wt{\ba} = B_{\ba}$ and therefore products along the $\wt{f}$ cycle with these elements are zero. It 
remains to adjust the product of $\wt{\alpha}$ and $\wt{\beta}$.

	 By 3.1.1 we have
	$p = c_{\ba}c_{\xi}c_{\alpha} A_{\alpha}$, and we see from the diagram that
and $A_{\alpha} = \alpha C \beta g(\beta)$ where $C$ is a monomial in the arrows of $Q$ of positive length and therefore, as an element of $R$, it belongs to 
the radical. 
We can replace the arrow
$\wt{\alpha}$ by 
$$\wt{\alpha}':= \wt{\alpha}(1-c_{\ba}c_{\xi}c_{\alpha}C)$$
and this has product zero with the arrow $\wt{\beta}$.

\bigskip
(a2) {\it Assume $p= \xi_{\alpha}= \alpha g(\alpha)f(g(\alpha))$ and $f(\alpha)$ is critical.}  Then we
use the diagram of 3.1.1 again, now taking
 $\alpha: i\to j$ and we set $\beta = f(g(\alpha))$ so that $p=\alpha \wt{\beta}$.  In this case, $\Gamma$ contains $i, k, x$ but not the vertex $j$.
 From this we see that the cycle of $\wt{f}$ containing $\wt{\alpha}$ is
 $$(\wt{\alpha} \ \beta \ \wt{f(\beta)} \ g(\beta)  \ \wt{\omega})
 $$
 here $\wt{f(\beta)}$  and $\wt{\omega}$ are socle elements and products with these along the cycle are zero, also after any adjustment. It remains to deal with
 the product of $\wt{\alpha}$ and $\beta$. 
 
  We have
 $$\xi_{\alpha} = c_{f(\alpha}c_{\omega}c_{\alpha}A_{\alpha}$$
  In this case, we see from the diagram that
 $A_{\alpha} = \wt{\alpha}\cdot C \cdot \beta$ where $C$ is a monomial of positive length. We set 
 $\wt{\alpha}' = \wt{\alpha}(1-c_{f(\alpha)}c_{\omega}c_{\alpha}C)$ and this can be taken as an arrow, and
 it satisfies  $\wt{\alpha}'\beta = 0$.

\bigskip

  Now consider  $p$ near a virtual arrow. \\
(a3) \ {\it Assume   $p=\zeta_{\alpha}$ so that $\wt{\alpha} = \alpha$,  and assume $\ba$ is virtual.}  We
have $\wt{\alpha} = \alpha$. Then $s(\alpha)$ and $s(f(\alpha))$ are in $\Gamma$ but $t(f(\alpha))$ is not in $\Gamma$.
In this case the virtual arrow $\ba$  cannot be a loop: 
Otherwise, using  part (3) of 3.1.2, we have  $\alpha: i\to j$ and both $i, j$ are in $\Gamma$. But then $f(\alpha): j \to i$ is an arrow of $\wt{Q}$ and
$f(\alpha)= \wt{f}(\wt{\alpha}) \neq f(\alpha)g(f(\alpha)).$

Now we use the diagram (1) of 3.1.2. We can assume that the virtual arrow $\ba$ is equal to $\xi$, that
is we take $\alpha: x\to k$. 
The set $\Gamma$ contains $x, k, i$ but does not contain $j$. Let $\beta: k\to i$.
We see that $\wt{f}$ has the cycle of length four, that is
$(\alpha \ \wt{f}(\wt{\alpha}) \ \beta \ \wt{\xi})$.
Moreover $\wt{\xi} = \xi \omega = B_{\xi}$ and it belongs to the socle. Therefore $\beta \wt{\xi} = 0$ and $\wt{\xi}\alpha=0$. 
The other two products need to be adjusted.
By 3.1.2 we have
$$\alpha \wt{f}(\wt{\alpha}) = c_{\xi}c_{\alpha} A_{\alpha} \ \mbox{ and } \ \wt{f}(\wt{\alpha})\beta = c_{\omega}c_{g(\alpha)} A_{g(\alpha)}.$$
Now, we can write $A_{\alpha} = \alpha\cdot C$ for a monomial $C$ of positive length between vertices in $\wt{Q}$,  and  luckily, we also have
$C\cdot \beta = A_{g(\alpha)}$, moreover $c_{\omega} = c_{\xi}$ and $c_{\alpha} = c_{g(\alpha)}$.  
Hence we can replace $\wt{f}(\wt{\alpha})$ by  $\wt{f}(\wt{\alpha})':= \wt{f}(\wt{\alpha})- c_{\xi}c_{\alpha} C$.

\bigskip

(a4) {\it  Assume $p=\xi_{\alpha}$, so that $\wt{\alpha} = \alpha g(\alpha)$, and assume $f(\alpha)$ is virtual. }
As in (a3), the virtual arrow cannot be a loop.
We use the diagram (1) of 3.1.2 
and we take $\alpha$ to be the arrow  $k\to j$. Then we have the following arrows of $\wt{Q}$ 
$$\wt{\alpha}: k\to i, \ \beta = \wt{f}(\wt{\alpha}): i\to x, \  \wt{\xi}: j\to j, \ \gamma: x\to k
$$
and they belong to the cycle of $\wt{f}$  of length four
$$(\wt{\alpha} \ \beta \ \wt{\xi} \ \gamma).
 $$
 Since $\wt{\xi} = \xi\omega = B_{\xi}$ is in the socle, the products with $\beta$ and $\gamma$ are zero. 
 We see from 3.1.2 that
 $$\wt{\alpha}\beta = c_{\omega}c_{\ba}A_{\ba}, \mbox{ and } \ \gamma \wt{\alpha} = c_{\xi}c_{\gamma}A_{\gamma}
 $$
 Moreover $A_{\gamma} = \gamma C$ and $C\beta = A_{\ba}$ and as well $c_{\gamma}= c_{\ba}$ and  $c_{\xi} = c_{\omega}$ . We replace $\wt{\alpha}$ by 
 $\wt{\alpha}' := \wt{\alpha} - c_{\omega}c_{\ba}C$, then the remaining products along the cycle of $\wt{f}$ are zero.

 \bigskip

(b)   {\it The case when $|p|=4$  and $p\neq 0$ in $\La$:} \ 
 Then by Lemma \ref{lem:7.5} we have $p = \wt{\alpha}\wt{\beta}$ where $\wt{\alpha} = \alpha g(\alpha)$ and $\wt{\beta} = \beta g(\beta)$ for
 $\beta = f(g(\alpha))$. 
 That is we can write $p= \xi_{\alpha} g(\beta) (= \alpha \zeta_{g(\alpha)})$, and we must have that $\xi_{\alpha}\neq 0$.
 This means that the arrow
 $f(\alpha)$ is critical or virtual. 
 
 (b1) {\it  Assume $f(\alpha)$ is a virtual loop.}  
 In this case we use the diagram (3) of 3.1.2, 
 with $\omega = f(\alpha)$ so that $\alpha: i\to j$ and $g(\alpha) (= f^2(\alpha)) : j\to i$. Then $\beta = \alpha$ and therefore $\wt{\alpha} = \wt{f}(\wt{\alpha})$ and it is a loop fixed by $\wt{f}$. We compute
 $$\wt{\alpha}^2 = c_{\alpha}c_{\omega}B_{\alpha}$$
 which is non-zero in the socle. This means that at $\wt{\alpha}$ we have a socle deformation.

 \medskip
 
 (b2) {\it Assume $f(\alpha)$ is virtual but not a loop.}  Then we use the diagram 3.1.2 with $\alpha: k\to j$, so that $\beta: i\to x$. Then 
 $\Gamma$ contains $k, i$ but does not contain $j, x$. We see that $\wt{f}$ has a  cycle of length two, namely
 $(\wt{\alpha} \ \wt{\beta})$. 
 Using the formulae in 3.1.2 we compute
 $$\wt{\alpha}\wt{\beta} = c_{f(\alpha)}c_{\alpha}B_{\alpha} = c_{f(\alpha)} c_{\alpha}\wt{\alpha} C$$
 where $C$ is a monomial of positive length from $i$ to $k$. Similarly
 $$ \wt{\beta}\wt{\alpha} = c_{f(\beta)}c_{\beta} B_{\beta} = c_{f(\beta)} c_{\bar{\beta}}B_{\bar{\beta}}
 $$
 using (3) of Definition \ref{def:3.3}. 
 Now $\bar{\beta}$ is in the $g$-orbit of $\alpha$ and we see $c_{\bar{\beta}} = c_{\alpha}$ and moreover
 $B_{\bar{\beta}} = C \wt{\beta}$. Furthermore $c_{f(\alpha)} = c_{f(\beta)}$
  We replace $\wt{\beta}$  by $\wt{\beta}':= \wt{\beta}-c_{\alpha}c_{f(\alpha)}C$.

 \bigskip
 
 (b3) {\it Assume $f(\alpha)$ is critical.}
 Then we use the diagram 3.1.1 with $\alpha: i\to j$, and $\beta= f(g(\alpha): k\to x$. Then $i,  k$ are in $\Gamma$ but $j, x$ are not in $\Gamma$. 
 The $\wt{f}$-cycle of $\wt{\alpha}$ is 
 $$(\wt{\alpha} \ \wt{\beta} \ \wt{\omega})$$
 and $\wt{\omega} = B_{\omega}$, hence the product of $\wt{\omega}$ with any arrow is zero.

 Using the calculations in 3.1.1 we have
 $$\xi_{\alpha}g(\beta) = c_{f(\alpha)}c_{\omega}c_{\alpha}B_{\alpha}$$
  We  factorise $B_{\alpha} = \wt{\alpha}C\wt{\beta}$ and $C$ is a monomial of positive length.  
  We can replace $\wt{\alpha}$ by $\wt{\alpha}':= \wt{\alpha}(1-c_{\alpha}c_{f(\alpha)}c_{\omega} C)$
  and the $\wt{\alpha}'\wt{\beta}=0$ (and $\wt{\omega}\wt{\alpha}'=0$).

(c) 
We determine now when the algebra $R=e\La e$ has only virtual  arrows, and then
verify that Condition (4) of Definition \ref{def:3.3} holds. 

(i) We show first that in this case, $R$ does not have a virtual arrow of type (b):

Suppose such an arrow $\wt{\alpha}$ say exists.  Then $m_{\alpha}\wt{n}_{\wt{\alpha}}=2$ and $\bar{\wt{\alpha}} = \wt{\ba}\in \wt{\cT}$. Then
$\wt{\ba} = \ba$. This must also be virtual and necessarily of type (b). Therefore also $\wt{\alpha}\in \wt{\cT}$ and then $\wt{\alpha}=\alpha$.

So we have two $f$-cycles of arrows in $Q$ which all remain arrows of $\wt{Q}$. 
If $\alpha, \ba$ are both loops then $\La$ must be local and $e=1$ which is excluded. So say $\alpha: i\to  j$ and $i\neq j$.
Then $i, j$ belong to $\Gamma$.Since $\wt{n}_{\wt{\alpha}}\leq 2$, the $g$-cycle of $\alpha$ cannot pass through any other vertex of $\Gamma$ and
$g(\alpha)$ is a path from $j$ to $i$. However $g^{-1}(\alpha) = f^{-1}(\ba)$ and it starts at some vertex in $\Gamma$. It follows that
$g(\alpha) = f^{-1}(\ba)$.

Assume (for a contradiction) that  $\ba$ is  a loop: \ 
Then $f(\ba): i\to j$ but $Q$ is 2-regular and then $f(\ba) = \alpha$. But then $f(\alpha)$ must be a loop at $j$ and
$Q$ has two vertices and moreover $Q=\wt{Q}$ and $e=1$ which is excluded. 

Then $Q$ has subquiver with three vertices which has arrows $\alpha, \ba, f(\ba)$ and $f^2(\ba)$. Now we can use the same reasoning for
$\ba$ and see that $f^{-1}(\alpha)$ which is $g^{-1}(\ba)$ is an arrow $k=t(\ba) \to i$. Then $f(\alpha): j\to k$ and $Q$ is the triangular quiver.
The algebra $\La$ has at least four virtual arrows and this is excluded in  4.3(2).

\bigskip

(ii) We have shown that if $R$ has only virtual  arrows then all 
arrows are virtual of  type (a), and hence they are
loops, and $R$ is local. 
Then $\wt{\alpha} = B_{\alpha}$ and $\wt{\ba} = B_{\ba}$. In particular $\wt{g}(\wt{\alpha}) = \wt{\alpha}$ and therefore $\wt{f} = (\wt{\alpha} \ \wt{\ba})$. We see that $R$ is the local algebra
as in \ref{ssec:4.1}(1) with both multiplicities equal to $1$. We also see that
condition (4) of \ref{def:3.3} holds.
\end{proof}

\bigskip

\subsection{The proof of 5.1 in the special cases}\label{ssec:5.1}
We consider the algebras which were excluded in the above proof.

\subsubsection{Idempotent algebras for a WSA as in 4.2 }
That is, $\La=KQ/I$ where the  quiver
$Q$  is of the form
\[
  \xymatrix{
%  \xymatrix@C=1pc{
    1
    \ar@(ld,ul)^{\alpha}[]
    \ar@<.5ex>[r]^{\beta}
    & 2
    \ar@<.5ex>[l]^{\gamma}
    \ar@(ru,dr)^{\sigma}[]
  }
%,
\]
and $f = (\alpha \ \beta \ \gamma)(\sigma)$ so that 
$g = (\alpha)(\beta \ \sigma \ \gamma)$. Let $m_{\alpha}=t\geq 2$ and $m_{\beta}=m$, and we can take $c_{\alpha} = \lambda$ and $c_{\beta}=1$.
By 4.2 (2a) (and Assumption \ref{ass:3.4}(1)), if $t=2$ then $m\geq 2$.
Furthermore, if $(t, m) = (3, 1)$ then $\lambda \neq 1$ (see 4.2(2b)).
There are two idempotent algebras $\neq \La$ to be considered, and we describe the result:\\
{\it
(1) Let $R=e_1\La e_1$. This gives a local algebra as in 4.1(1). In particular for $m=1$ we have $R\cong K[x]/(x^t)$.\\
(2) Let $R=e_2\La e_2$, then $\wt{\cT} = \{ \sigma\}$ and we get the algebras as in 4.1(2b). When $t>2$ it is of semidihedral type, 
and if $t=2$ it is a socle deformation of an algebra of semidihedral type. }

We omit details for (1), but we give details for (2), to show how a socle deformation occurs.
Hence let  $R:= e_2\La e_2$. This algebra has
quiver
\[
 \xymatrix{
         2 \ar@(dl,ul)[]^{\wt{\gamma}} 
    \ar@(ur,dr)[]^{\sigma}
 }
\]
where
$\wt{\gamma} = \gamma\beta$.
The permutations are
$$\wt{f} = (\wt{\gamma})(\sigma)  \ \mbox{ and } \ \wt{g} = (\wt{\gamma} \ \sigma).
$$
In this case  we have $\wt{\cT} = \{ \sigma\}$. We write down the type (1)
relations of Definition \ref{def:3.3}. The first one is
$$\sigma^2 = A_{\wt{\gamma}} = (\wt{\gamma}\sigma)^{m-1}\wt{\gamma}.
$$
Next, $\wt{\gamma}^2 = \gamma\beta\gamma\beta =0$ provided $t>2$, by the
zero relations for $\La$. Assume now $t=2$, then using the relations for
$\La$ we see
$$\gamma\beta\gamma\beta = \lambda \gamma\alpha\beta = \lambda \gamma A_{\beta}
= \lambda B_{\gamma}
$$
which is non-zero and spans the socle of $R$. That is, we get an algebra
as in 4.1(2b) when $t>2$. If $t=2$ we get a socle deformation of
such an algebra.

\bigskip

\subsubsection{Idempotent algebras when $\La$ is a WSA as in 4.3(3)}

Then the quiver is triangular, and 
we have   $m_{\bullet} = (2, 2, 1)$ and $c_{\bullet} = (\lambda, 1, 1)$.  The arrows $\alpha_3, \beta_3$
are virtual, and up to labelling we have to consider four idempotent algebra. 
We describe the result, the details are straightforward and are omitted.\\
{\it  (1) If  $e= e_1+e_2$ and $R=e\La e$ then $R$ is a Brauer graph algebra with one virtual loop.\\
(2) If $e=e_1 + e_3$ then again $R=e\La e$ is a Brauer graph algebra. In this case, the virtual arrows of $\La$ are not virtual
as arrows of $R$,\\
(3) The algebra $e_1\La e_1$ is a local hybrid algebra as in 4.1(1).\\
(4) The algebra $e_2\La e_2$ is a 4-dimensional algebra of dihedral type, as in 4.1(1).}
\bigskip

\begin{remark}\label{rem:5.5} \
(a)  Suppose $\wt{\alpha}$ is an arrow of $\wt{Q}$ starting at
        $i$. We must show
        that $\wt{m}_{\wt{\alpha}}\wt{n}_{\wt{\alpha}}=1$ only occurs when the vertex
        is biserial and $\wt{\alpha}$ is a loop.

       We have $\wt{n}_{\wt{\alpha}}=1$ if and only if
        $\wt{\alpha}$ is the product of all arrows in the $g$-cycle of $\alpha$, hence is a loop. If in addition $\wt{m}_{\wt{\alpha}}=1$
        then $\wt{\alpha} = B_{\alpha}$ and clearly $\wt{\alpha}\wt{\ba}=0$
        and $\wt{\ba}\wt{\wt{\alpha}}=0$. To see that $i$ is biserial we need
	$\wt{\ba}$ is not in $\wt{\cT}$. This is clear if $\alpha$ 
	is a loop since then $\alpha = B_{\alpha}$ and $\alpha$ is virtual 
	of type (a). Suppose $\alpha$ is not a loop and
	$\wt{\alpha} = B_{\alpha}$. The last arrow in $B_{\alpha}$ is $f^2(\ba)$ and 
	it does not start at a vertex of $\Gamma$ and therefore
	$\wt{\ba}$ cannot be in $\wt{\cT}$.

        (b) \ The algebra $e\La e$ is symmetric, therefore the exceptions in 
	Assumption \ref{ass:3.4} cannot occur.
\end{remark}

\bigskip

We will now show that every hybrid algebra, such that $\cT \neq Q_1$, 
occurs as an
idempotent algebra of some weighted surface algebra.
This generalizes the main result of \cite{BGA} where this was done
for the case of Brauer graph algebras. 
As in \cite{BGA}, our tool is the $*$-construction which we will now introduce.

\subsection{The $*$-construction}\label{ss:5.2}
Let $H$ be a hybrid algebra such that 
$\cT \neq Q_1$,  
say $H= H_{\cT}(Q, f, m_{\bullet}, c_{\bullet})$, and let $g$ be the permutation associated to $f$.
The $*$-construction gives  a triangulation quiver $(Q^*, f^*)$ which contains $Q_0$, and furthermore, contains all arrows in $\cT$. 

The idea is to keep the arrows of $\cT$ as they are, but 
split each arrow which is not in $\cT$, and add extra arrows
in order to create  triangles. 
With this, one has weighted surface algebras with $m^*, c^*$ extending $m, c$. Explicitly, define
$$Q_0^*:= Q_0\cup \{ x_{\alpha}\}_{\alpha \in Q_1\setminus \cT}, \ \ 
Q_1^*:= \cT\cup \{ \alpha', \ \alpha'', \ \ve_{\alpha}\}_{\alpha\in Q_1\setminus \cT}$$
For $\beta\in \cT$ we set $s^*(\beta)=s(\beta)$ and $t^*(\beta)=t(\beta)$.
Let $\alpha$ be an arrow which is not in $\cT$. Then we set
$$\begin{aligned} 
        s^*(\alpha'):= s(\alpha), & \ &  t^*(\alpha'):= x_{\alpha}, & \ &
        s^*(\alpha''):= x_{\alpha}, & \ & t^*(\alpha)'':= t(\alpha)
        \cr
 s^*(\ve_{\alpha}) = x_{f(\alpha)},& \  & t^*(\ve_{\alpha}) = x_{\alpha}.
 \end{aligned}
$$
Next we  define the permutation $f^*$ on
 $Q^*$. If $\beta\in \cT$ then we take $f^*(\beta) = f(\beta)$, and  define
$$f^*(\alpha''):= f(\alpha)', \ \ f^*(f(\alpha)'):= \ve_{\alpha}, \ \ f^*(\ve_{\alpha}) := \alpha''.
$$
Then $(Q^*, f^*)$ is a triangulation quiver.

This determines the permutation $g^*$, explicitly it is as follows.
First, if the arrow $\alpha$ of $Q$ is not in $\cT$ then $g^*(\alpha') = \alpha''$.
The arrows starting at $t(\alpha'')$ in $Q$ are $f(\alpha)$ and $g(\alpha)$, and $g^*(\alpha'')$ depends on whether or not
$g(\alpha)$ is in $\cT$, that is
$$g^*(\alpha'')  \ = \left\{\begin{array}{ll} g(\alpha)' & g(\alpha)\not\in \cT,\cr
g(\alpha) & \mbox{ else}.
\end{array}
\right.
$$
Finally, $g^*(\ve_{\alpha}) = \ve_{f^{-1}(\alpha)}$ for any $\alpha\in Q_1\setminus \cT$.
The cycles of $g^*$ are obtained from the cycles of $g$ by replacing each $\alpha$ in $Q_1\setminus \cT$ by $\alpha', \alpha''$, together with cycles only
containing arrows of the form $\ve_{\alpha}$.  On the cycles without $\ve$-arrows, we take the same multiplicity function and parameter function as for $H$.
On the $\ve$-cycles we may choose multiplicities and parameters arbitrarily. We take them equal to $1$ unless
when some arrow $\ve_{\gamma}$ is required to be not virtual or critical, then we choose $m_{\ve_{\gamma}} \geq 3$, or when
some non-zero scalar factor needs to be specified, we may choose $c_{\ve_{\gamma}}$ differently, depending  on the context.
This defines then a weighted surface algebra $\Lambda = \La(Q^*, f^*, m^*, c^*)$.
In fact, this is a choice, we could equally well  apply the $*$ construction also to triangles in $\cT$.

Note that when $\cT = Q_1$, the construction does not do anything, and $H$ is already a weighted surface algebra as in \cite{WSA-GV}. The case when 
$H$ is local and $\cT = Q_1$ is discussed in 4.1(2a), and this
is not a weighted surface algebra by the definition in \cite{WSA-GV}.

\bigskip

\begin{example}\label{ex:5.7}
We illustrate the $*$-construction.\\
(1)
A  loop $\alpha$ in $Q$ fixed by $f$ which does not belong to $\cT$ is replaced
in $Q^*$ by the subquiver
\[
  \xymatrix{ 
     x_{\alpha} \ar@(dl,ul)[]^{\varepsilon_{\alpha}} \ar@/^1.5ex/[r]^{\alpha''}
       & s(\alpha) \ar@/^1.5ex/[l]^{\alpha'}
  }
\]
which is an orbit of  $f^*$.
\smallskip

(2)
An $f$-cycle in $Q$ which does not belong to $\cT$ of the form
\[
  \xymatrix{ 
     a \ar@(dl,ul)[]^{\alpha} \ar@/^1.5ex/[r]^{\beta} & b \ar@/^1.5ex/[l]^{\gamma}
  }
\]
is replaced in $Q^*$ by the quiver
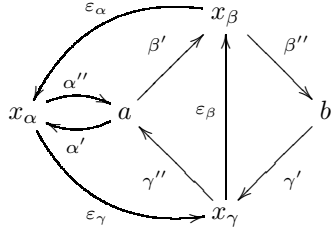
\[
  \xymatrix{ 
     && x_{\beta} \ar[rd]^{\beta''} \ar@/_4ex/[lld]_{\varepsilon_{\alpha}} \\
     \!\! x_{\alpha} \ar@/^1.5ex/[r]^{\alpha''}  \ar@/_4ex/[rrd]_{\varepsilon_{\gamma}}
       & a \ar@/^1.5ex/[l]^{\alpha'} \ar[ru]^{\beta'} && b \ar[ld]^{\gamma'} \\
     && x_{\gamma} \ar[lu]^{\gamma''} \ar[uu]^{\varepsilon_{\beta}}      
  }
\]
with $f^*$-orbits
$(\alpha'' \ \beta' \ \varepsilon_{\alpha})$,
$(\gamma'' \ \alpha' \ \varepsilon_{\gamma})$ and $(\beta'' \ \gamma' \ \varepsilon_{\beta})$.
\smallskip

(3) \
Suppose $f$ has a 4-cycle
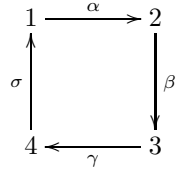
\[
%  \xymatrix{
  \xymatrix@C=3pc@R=3pc{
    1
      \ar[r]^{\alpha}
  &
    2
      \ar[d]^{\beta}
  \\ 
    4
      \ar[u]^{\sigma}
  &
    3
      \ar[l]^{\gamma}
  } 
\]
Then the corresponding part of $Q^*$ is of the form
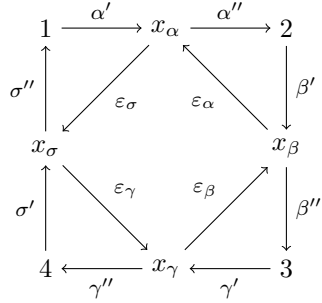
\[
\begin{tikzpicture}
%[->,scale=1]
[->,scale=.8]
%[->,scale=.8]
\coordinate (1) at (-2,2) ;
\coordinate (2) at (2,2) ;
\coordinate (3) at (2,-2) ;
\coordinate (4) at (-2,-2) ;
\coordinate (up) at (0,2) ;
\coordinate (down) at (0,-2) ;
\coordinate (right) at (2,0) ;
\coordinate (left) at (-2,0) ;
%
%\fill[rounded corners=3mm,fill=gray!20] (1) -- (up) -- (left) -- cycle;
%\fill[rounded corners=3mm,fill=gray!20] (2) -- (right) -- (up) -- cycle;
%\fill[rounded corners=3mm,fill=gray!20] (3) -- (down) -- (right) -- cycle;
%\fill[rounded corners=3mm,fill=gray!20] (4) -- (left) -- (down) -- cycle;
%
\node (1) at (-2,2) {$1$};
  \node (2) at (2,2) {$2$};
\node (3) at (2,-2) {$3$};
\node (4) at (-2,-2) {$4$};
\node (up) at (0,2) {$x_{\alpha}$};
\node (down) at (0,-2) {$x_{\gamma}$};
\node (right) at (2,0) {$x_{\beta}$};
\node (left) at (-2,0) {$x_{\sigma}$};
\draw
(1) edge node[above]{\footnotesize$\alpha'$} (up)
(up) edge node[above]{\footnotesize$\alpha''$} (2)
(2) edge node[right]{\footnotesize$\beta'$} (right)
(right) edge node[right]{\footnotesize$\beta''$} (3)
(3) edge node[below]{\footnotesize$\gamma'$} (down)
(down) edge node[below]{\footnotesize$\gamma''$} (4)
(4) edge node[left]{\footnotesize$\sigma'$} (left)
(left) edge node[left]{\footnotesize$\sigma''$} (1)
(right) edge node[below left]{\footnotesize$\varepsilon_{\alpha}$} (up)
(down) edge node[above left]{\footnotesize$\varepsilon_{\beta}$} (right)
(left) edge node[above right]{\footnotesize$\varepsilon_{\gamma}$} (down)
(up) edge node[below right]{\footnotesize$\varepsilon_{\sigma}$} (left)
;
\end{tikzpicture}
\]
\end{example}

\bigskip

\begin{theorem} \label{thm:5.8} Assume $H$ is a hybrid algebra, such that $\cT\neq Q_1$.  
Then there is a weighted surface algebra  $\La$ and an idempotent $e$ of $\La$ such that
$H$ is isomorphic to a block component of $e\La e$.
\end{theorem}

\bigskip
{\it Proof }  Given $H= H_{\cT}(Q, f, m_{\bullet}, c_{\bullet})$. 
We let $(Q^*, f^*)$ and $\La$ as constructed above.
Now let $e$ be the idempotent $e:= \sum_{i\in Q_0} e_i$. We want to show that $e\La e$ is isomorphic to $H$. \\

We have three algebras, the given algebra is $H= KQ/I$, next we
have the weighted surface algebra 
$\La = KQ^*/I^*$ associated to the triangulation quiver $(Q^*, f^*)$ as 
introduced above. Furthermore,  we have the idempotent algebra $e\La e$. 
By  Theorem \ref{thm:5.1} we know  that
it
has a presentation $K\wt{Q^*}/\wt{I^*}$ and that it is a hybrid algebra. 

Since $e=\sum_{i\in Q_0} e_i$, the quiver $\wt{Q^*}$ has vertices $(\wt{Q^*})_0 = Q_0$.
The arrows of 
$\wt{Q^*}$ are obtained by contracting paths of $Q^*$ 
of shortest length between vertices in $Q_0$. The arrows of $Q^*$ are
\begin{enumerate} [(1)] 
	\item \emph{the arrows of $\cT$,}
	\item \emph{ arrows $\alpha', \  \alpha''$ and $\ve_{\alpha}$ for each 
		arrow $\alpha \in Q_1\setminus  \cT$.}
\end{enumerate}

The arrows of $Q^*$ starting at some vertex in $Q_0$ are therefore
the $\alpha$ in $\cT$, and the $\alpha'$ when $\alpha\not\in \cT$. 
If $\alpha \in \cT$ then $\wt{\alpha}=\alpha$, and if $\alpha\not\in \cT$ then 
$\wt{\alpha'} = \alpha'\alpha''$. 
So $\wt{Q^*}_1$ is the set of $\wt{\alpha}$ for $\alpha \in \cT$ and $\wt{\alpha'}$ for $\alpha \in Q_1\setminus \cT$.

The set of triangles $\wt{\cT}$ of the algebra $K\wt{Q^*}/\wt{I^*}$ 
 consists therefore  of the set
$\{ \wt{\alpha} \mid \alpha \in \cT\}$
(see part (c) in the proof of Theorem  \ref{thm:5.1}).
We define a surjective algebra map $\psi: K\wt{Q^*} \to H$ by
$\psi(e_i) = e_i$ and
if $\wt{\gamma}$ is an arrow of $\wt{Q^*}$ then 
$$\psi(\wt{\gamma}) = \left\{\begin{array}{ll} \gamma, & \mbox{if} \ \wt{\gamma} =\gamma\cr
	\alpha,  & \mbox{if} \ \wt{\gamma} = \alpha', 
\end{array}\right.
$$
and extending to products and linear combinations.

We show now that $\psi(\wt{I})=0$ (that is $\psi$ induces an algebra homomorphism from $e\La e$ to $H$).
First we observe that  $\psi$ takes any  submonomials of $B_{\wt{\gamma}}$ starting and ending at vertices
in $Q_0$ to its  'contraction', replacing each subpath of the form $\alpha'\alpha''$ by $\alpha$, and leaving each
$\gamma\in \cT$ unchanged.

\bigskip

(a) We consider relation (1) of Definition \ref{def:3.3}.  Assume $\wt{\gamma} \in \wt{\cT}$, then   we have
$\wt{\gamma} \wt{f(\gamma)} = \wt{\gamma}\wt{f(\gamma)} = c_{\bar{\wt{\gamma}}} A_{\bar{\wt{\gamma}}}$.
By the above observation we see see that
$\psi$ preserves this identity.
Now consider  an arrow of the form  $\wt{\alpha'}$  for $\alpha\in Q_1$ and not in $\cT$. Then we have
$$\wt{\alpha'}\wt{f}(\wt{\alpha'})  = \alpha'\alpha''\cdot f(\alpha)'f(\alpha)''
\leqno{(*)}$$
    By definition,
$\psi(\wt{\alpha'})\psi(\wt{f(\alpha)'}) = \alpha f(\alpha) = 0.$
By our convention, we can make sure that $\ve_{f(\alpha)} (= f(\alpha'))$ is not virtual or critical
Then the path $\alpha'\alpha''\cdot f(\alpha)'f(\alpha)''$ is zero, by Lemma \ref{lem:7.1} (see Appendix).

(b) Next consider a loop of the form $\wt{\alpha'}$ for $\alpha \in Q_1$ and $\wt{\alpha'}$ not in $\wt{\cT}$, with
$\wt{\alpha'} = \wt{f}(\wt{\alpha'})$. 
Then we have $f(\alpha) = \alpha$ and $\alpha^2=0$. Now 
$$\wt{\alpha'}^2 = \alpha'\alpha''\cdot \alpha'\alpha''. \leqno{(*)}
$$
By definition
$\psi(\wt{\alpha'})^2  = \alpha^2 = 0$.
The subquiver of $Q^*$ constructed from  a loop $\alpha$ fixed by $f$ is shown in Example \ref{ex:5.7}(1).
We have
$$\alpha'\alpha''\alpha'\alpha'' = c\alpha'A_{\ve_{\alpha}}\alpha'' \leqno{(\dagger)}
$$
where $c=c_{\ve_{\alpha}}\neq 0$.
We may choose $c$ and we may also choose $m_{\ve_{\alpha}}$. We take $m_{\ve_{\alpha}}$ large enough so that $\ve_{\alpha}$ is
not virtual or critical, and then
$(\dagger)$ is zero.

(c) Now consider the relations   (2) and (2') of Definition \ref{def:3.3} when $\alpha$ and $\ba$ (respectively $g(\alpha)$) are in $\cT$. Then also $f(\alpha)$ is in $\cT$
and this part of the quiver,  the map $\psi$ is an identification, so the relations are preserved.
Otherwise the elements are mapped to zero by (1) of Definition \ref{def:3.3} .
The socle relations (3) follow automatically.

To complete the proof it suffices to establish that
$e\La e$ and $H$ have the same dimensions.
For any vertex $i$, 
the dimension of $e_iH$ is $m_{\alpha}n_{\alpha} + m_{\ba}n_{\ba}$, and
it is the same as that of $e_i(e\La e)$.
$\Box$

\begin{example}\label{ex:5.9}
Let $\La$ be the local algebra with arrows $\alpha, \beta$ and
$$f=(\alpha)(\beta), \ \ g = (\alpha \ \beta).
$$
We take $\cT=  \{ \beta \}$ with $m_{\bullet}=1$ and $c_{\bullet}=c$ so that $\alpha$ is
        virtual.
The relations are
        $$\beta^2 = cA_{\alpha}, \ \ \alpha^2 = 0
$$
and the zero relation $\alpha\beta\alpha = 0$.
We apply the $*$ construction to $\alpha$. This gives
        the algebra $\La^*$ with quiver
\[
% \xymatrix{
%  \bullet \ar@(dl,ul)[]^{\alpha} \ar@<+.5ex>[r]^{\sigma}
%     \save[] +<0pc,3mm> *{1} \restore
%   & \bullet \ar@<+.5ex>[l]^{\gamma} \ar@(ur,dr)[]^{\beta}
%     \save[] +<0pc,3mm> *{2} \restore
% }
 \xymatrix{
         x_{\alpha} \ar@(dl,ul)[]^{\ve_{\alpha}} \ar@<+.5ex>[r]^{\alpha''}
   & i \ar@<+.5ex>[l]^{\alpha'} \ar@(ur,dr)[]^{\beta}
 }
\]
Take $m_{\ve_{\alpha}}=4$.
         We may write down the relations defining $\La$, for simplicity write
        $\ve = \ve_{\alpha}$ and $d= c_{\ve}$.

$$\begin{aligned}
        \alpha''\alpha'=& dA_{\ve}, &\ & \ \alpha'\ve = cA_{\beta}, &\ & \ve\alpha''=& cA_{\alpha''}, &\ & \beta^2 =& cA_{\alpha'}
\end{aligned}
$$
        together with the zero relations, in particular $\alpha'\alpha''\alpha'=0$.

Now consider the idempotent algebra $e\La e$, we want this to be isomorphic
to $H$.
        By Theorem \ref{thm:5.1} it has a presentation $K\wt{Q}/\wt{I}$ where
        $\wt{Q}$ is the quiver with two loops $\wt{\alpha'}$ and $\wt{\beta}$, and $\wt{\alpha'} = \alpha'\alpha''$, and $\wt{\beta} = \beta$. This has
        relations
        $$\wt{\beta}^2 = cA_{\wt{\alpha'}}, \ \ 
        (\wt{\alpha'})^2 = 0
$$
\end{example}

\bigskip

\begin{remark}\label{rmk:5.10}
The  algebra 
 $\La$ in
the proof of Theorem \ref{thm:5.8} is a WSA and hence is symmetric, so
	it is not one of the exceptions in Assumption \ref{ass:3.4}.
\end{remark}

\begin{lemma}\label{lem:5.11} Assume $H$ is a hybrid algebra. Then $H$ is tame and symmetric.
\end{lemma}

{\it Proof} \
We have proved that  any hybrid algebra is an idempotent algebra
of a (general) weighted surface algebra. Weighted surface algebras are tame and symmetric
(see \cite{WSA-GV}), and it is well known that idempotent algebras of tame symmetric algebras are tame and symmetric.
$\Box$

\bigskip

\section{Stable Auslander-Reiten  components }

This section is more general, here we assume  $\La$ is a tame symmetric algebra such that its Gabriel quiver is 2-regular. 
We can take $\La$ to be basic, with an admissible presentation 
$\La = KQ/I$ and hence  $Q$ is 2-regular. 

For background we refer to Chapter 4 in \cite{Benson}.

The Auslander-Reiten (AR) quiver $\Gamma_{\La}$ of an algebra $\La$ is the graph where the vertices correspond to isomorphism types
of indecomposable $\La$-modules, and where the arrows are labelled in terms of irreducible maps.
For our context it is most relevant that this quiver encodes Auslander-Reiten (AR) sequences, also known as almost split sequences.

A short exact sequence $0\to M\to E \stackrel{\sigma}\to N\to 0$ is an  AR sequence if $M$ and $N$ are indecomposable, the map $\sigma$ does not split, and moreover given any module $N'$ and a map $\rho:  N'\to N$ which is not a split epimorphism, then $\rho = \psi\circ \sigma$ for some $\psi: N'\to E$.
It was proved by Auslander and Reiten \cite{AR} that for any indecomposable non-projective module $N$, such a sequence exists, and it is unique up to
isomorphism of short exact sequences. The module $M$ is denoted by $\tau(N)$ and $\tau$ is known as  Auslander-Reiten translation.
 The arrows in $\Gamma_{\La}$ are then as follows: For $N$ indecomposable non-projective,
the number of arrows $X\to N$ is the multiplicity of $X$ as a direct summand of $E$
(which usually is $\leq 1$.
For $M$ indecomposable and not injective, there is an almost split sequence starting with $M$. Then the number of arrows from $M$ to $X$
is the multiplicity of $X$ as a direct summand of $E$.

We assume the algebra is symmetric, so that projectives and injectives are the same.
In this case we have $\tau \cong \Omega^2$.
The only almost split sequence in which an indecomposable projective $P_i$ corresponding to the simple module $S_i$ can occur, is what we call standard sequence
$$0\to \Omega(S_i) \to P_i\oplus {\rm rad}(P_i)/S_i \to \Omega^{-1}(S_i)\to 0$$

We assume that $\La$ is symmetric, then the stable AR-quiver $_s\Gamma_{\La}$ is obtained from $\Gamma$ by removing the vertices corresponding to the indecomposable
projective modules.
The stable AR quiver is a translation quiver, where $\Omega^2$ acts as translation. 
The graph structure of a component of $_s\Gamma_{\La}$ is described by
Riedtmann's structure theorem.

For each component $\cC$ if $\Gamma_{\La}$, its stable part is a component of
$_s\Gamma_{\La}$, and for $\La$ of infinite type, $\cC$ is either a stable
tube $\cC\cong \bZ A_{\infty}/(\tau^r)$ (if it contains a periodic module
\cite{HPR}), or it is an (acyclic) quiver of the form $\cC = \bZ\Delta$.

The main  tool to identify the graph structure
of $\cC$  are subadditive functions, by  applying  the
classification theorem of \cite{HPR}. For the case of group algebras of 
finite groups, this was done by Webb \cite{W}, and Okuyama presented a new
approach \cite{O}. 
We use the version from Section 3 of \cite{ES} where
this 
is generalized to selfinjective algebras. The identification method is
then described as follows.

We say that $\La$ has enough
periodic modules if
	for each indecomposable non-projective $M$ there is a module
	$W$ with $W\cong \tau(W)$, such that $\ul{\rm Hom}_{\La}(W, M)$ is non-zero.
	Here $\ul{\rm Hom}_{\La}(X, Y) = {\rm Hom}_{\La}(X, Y)/P(X, Y)$ where $P(X, Y)$ is the 
	subspace of maps which factor through some projective module.
	 Note that $\tau$-periodic is the same as $\Omega$-periodic for symmetric algebras.

\begin{proposition}\label{prop:6.1} Assume $\La$ has enough periodic modules. 
Let $\Theta$ be the stable component containing some
	indecomposable non-projective module $M$, let $W$ be
	as above. Then $d_W:= \dim \ul{\rm Hom}(W, -)$ defines
	an additive function on $\Theta$, hence $T$ is either
	Dynkin or Euclidean or one of the infinite trees $A_{\infty}, A_{\infty}^{\infty}, D_{\infty}$.
\end{proposition}

When $\Theta$ contains a periodic module then $T\cong A_{\infty}$ (for $\La$ of infinite type), see \cite{HPR}. If $\Theta$ contains
no periodic modules then both $M$ and its syzygy $\Omega(M)$ are not summands
of $W$, and then $d_W$ is an additive function, by \cite[Lemma 3.2]{ES}. 
The problem is how to find such module $W$ when 
modules in  
$\Theta$ are not periodic.

\subsection{Finding modules $W$}

Assume $\La$ is tame and symmetric. 
Furthermore, we assume that the Gabriel quiver of $\La$ is 2-regular. This means that every component
$\cS$ of the separated quiver is of the form $\wt{A}_n$ for some $n$. 

We recall the definition of  the separated quiver of an algebra.
If $Q$ is the quiver of the algebra and 
has vertices labelled by $1, 2, \ldots, r$ then the separated quiver
$Q_s$ has vertices $\{ 1, 2, \ldots, r, 1', 2', \ldots, r'\}$. The arrows
of $Q_s$ are given by
$\alpha: i\to j'$  whenever $\alpha: i\to j$ is an arrow in $Q$.
If $Q$ is a 2-regular quiver then there are two arrows starting at
each of $1, 2, \ldots, r$ of $Q_s$, and there are two arrows ending at each of
$1', 2', \ldots, r'$ of $Q_s$. Hence each component of $Q_s$ is isomorphic to
$\wt{A_n}$ for some $n$ (possibly a Kronecker quiver).

By the well-known classification of indecomposables
of such a quiver, there is 
a 1-parameter family of $K\cS$-modules $W_{\lambda}$ (for $\lambda\in K^*$) of $\tau$-period 1, all of dimension equal to the number of vertices of $\cS$.
Note that they have  radical length two.

\bigskip

The modules $W_{\lambda}$ can be viewed as $\La$-modules (by letting 
the square of the radical act as zero). By \cite{CB} 
they must be (almost all) periodic as $\La$-modules since the
algebra is tame, still of $\tau$-period 1, and therefore of   $\Omega$-period $2$ for $\La$. 
The same holds for an arbitrary component of the separated quiver.  There is some $\lambda\in K^*$ such that the $W_{\lambda}$ for each component
are periodic of period 2 as modules for $\Lambda$. 
Define
$$W_0:= \oplus_{\cS} W_{\lambda, \cS} \ \ \mbox{ and} \ W:= W_0\oplus \Omega_{\La}(W_0)
$$
Then $W$ is a periodic $\La-$module with $\Omega(W)\cong W$.

We take this module $W$, and let $d_W$ as above. By construction, 
$W_0$ has radical length $=2$ and ${\rm soc}(W_0) \cong W_0/{\rm rad} W_0 \cong \oplus_{i\in Q_0} S_i$.
We may take a set of minimal generators $\{ v_1, \ldots, v_n\}$ of $W_0$ such that $v_i = v_ie_i$.
Then we can take a basis of ${\rm soc } (W_0)$, of the form $w_1, \ldots, w_n$ such that $w_i= w_ie_i$.
Then if for some $i$ the arrows in $Q$ starting at $i$ are
$\alpha, \ba$ ending at $j, k$ then $v_i\alpha$ and $v_i\ba$ are non-zero, and are scalar multiples of $w_j, w_k$ respectively
(and we may have $j=k$).

\begin{lemma} \label{lem:6.2} Assume $M$ is indecomposable and not projective, and $\ul{\rm Hom}(W_0, M)=0$. Then
${\rm Hom}(W_0, M) \cong {\rm soc}(M)$.
\end{lemma}

\medskip

\begin{proof} \  (a)  
We define a homomorphism $\phi: {\rm soc}(M) \to {\rm Hom}_{\La}(W_0, M)$. 
We fix a $K$-basis for  ${\rm soc}(M)$ of the form 
$\{ m_{i,\nu(i)} \mid i\in Q_0, 1\leq \nu(i) \leq t_i\}$ where
$m_{i, \nu(i)} = m_{i, \nu(i)}e_i$.  Now define a linear
map
$$f_{i\nu(i)}: W_0 \to M$$
by $f_{i\nu(i)}(v_j) = \delta_{ij}m_{i\nu(i)}$ and $f_{i\nu(i)}(w_x)=0$. 
This defines a $\La$-module homomorphism. Now define $\phi(m_{i\nu(i)} = 
f_{i\nu(i)}$.

(b) We show that $\phi$ is injective:
Suppose $\phi(m)=0$ where $m=\sum_{i, \nu(i)}c_{i\nu(i)}m_{i\nu(i)}$ with $c_{i\nu(i)} \in K$, so $\phi(m) = \sum_{i, \nu} c_{i\nu(i)}f_{i\nu(i)}$.
Applying  this to some generator of $W_0$ gives 
$$0 = \phi(m)(v_j) = \sum_{\nu} c_{j\nu(j)}m_{j\nu(j)} 
$$
and since the $m_{j\nu(j)}$ are linearly independent it follows that all
$c_{j\nu(j)}$ are zero. Hence $m=0$.

(c) We show that $\phi$ is surjective.
Suppose there is some homomorphism  $f: W_0\to M$. 
It suffices to show that $f({\rm soc}(W_0))=0$:  if so 
then $f$ factors through $W_0/{\rm soc}(W_0)$ which is semisimple, and the image is contained in the socle. Then $f(v_j) = \sum_{\nu} c_{j\nu(j)} m_{j\nu(j)}$ with $c_{j\nu(j)}\in K$ for each
$j$ and  $f = \sum_{i \nu(i)} c_{i\nu(i)} f_{i\nu(i)}$, which  is in the image of 
$\phi$.

Assume false, then we may assume $f(w_r)$ is non-zero for some $r$.
We consider the diagram
$$\CD  0 @>>> W_0 @>{\iota}>> \La @>>> \Omega(W_0) @>>> 0 \cr
&& @V{f}VV \cr
&& M 
\endCD
$$
where $\iota$ is  the inclusion map.
Since $\ul{\rm Hom}(W_0, M)=0$,  it follows that $f$ must factor through $\iota$, so there is $h: \La \to M$ such that
$$f = h\circ \iota.
$$
Now,  $\iota(w_r)$ must span the socle of the copy of $e_r\La$ of $\La$
and we have $f(w) = h(\iota(w)) \neq 0$. Therefore the restriction of
$h$ to $e_r\La$ is non-zero, and then it is a split monomorphism, since $e_r\La$ is also injective.
This is not possible since $M$ is indecomposable and not projective.
So we have a contradiction.
\end{proof}

\bigskip

For the next part we will use an explicit injective hull of $W_0$. Note that its socle is multiplicity-free, and that every simple module occurs. 
We know that $W_0\cong \Omega^2(W_0)$, hence there is an exact sequence
$$0\to W_0 \to \La \to \Omega(W_0)\to 0
$$
and moreover since $W_0$ has radical length $= 2$, it is contained in the second socle of $\La$.

\bigskip

\begin{lemma} \label{lem:6.3} Assume $M$ is indecomposable and not projective , such that
$\ul{\rm Hom}(\Omega (W_0), M)=0$. Then ${\rm Hom}(\Omega(W_0), M) \cong {\rm rad}(M)$.
\end{lemma}

\bigskip

\begin{proof} \ \ 
(a) We show first that every $f: \Omega (W_0) \to M$ maps into 
	the radical of $M$. 
Suppose  there is some $f$ and  $f(x)$ is not in ${\rm rad}(M)$ for some 
	$x\in \Omega(W_0)$, then we may assume $f(x) = f(x)e_i$.
Since $f$ is zero in $\ul{\rm Hom}(\Omega (W_0), M)$, there is $h: \La \to M$ and $f = h\circ \iota$. In particular there is
$z = ze_i \in \La$ and $h(z) = f(x)$. Then $z$ must be a generator of $\La$ and $z\La \cong e_i\La$. The restriction of
$h$ to $z\La$ must split since $e_i\La$ is projective, and $M$ has a projective direct summand,  a contradiction.

We identify ${\rm Hom}(\Omega(W_0), M)$ with the set of 
$f: \La \to {\rm rad}(M)$ which take $W_0$ to zero. 

(b) \ We claim that if $f$ maps into the radical of $M$ then $f({\rm soc}_2(\La) =0$, and hence $f(W_0)=0$.
Let $f(e_i) = m = me_i$ in the radical of $M$. Then we can write $m = z\beta + z^*\beta^*$ where $\beta, \beta^*$ are
the arrows of $Q$ ending at $i$, and where $z$ and $z^*$ are elements of $M$.

Suppose there is some element $A$ in ${\rm soc}_2(\La)$ with $mA\neq 0$, say  $z\beta A\neq 0$. 
Then in particular $\beta A$ is non-zero in the socle of $e_j\La$ (for $j=s(\beta)$). 
 It follows that the submodule $z\La$ of $M$ is isomorphic to $e_j\La$. But $e_j\La$ is injective, and hence
 is a direct summand of $M$. This is a contradiction since $M$ is assumed to be indecomposable and not projective
 (hence injective).

(c) We define a homomorphism $\phi: {\rm rad}(M) \to {\rm Hom}_{\La}(\Omega(W_0), M)$, as in the proof of Lemma \ref{lem:6.2}. 

Take a basis of ${\rm rad}(M)$ of the form 
	$\{ m_{i\nu(i)} \mid i\in Q_0,  1\leq \nu(i) \leq s_i\}$ with 
	$m_{i\nu(i)} \in Me_i$. Then define on the generators of $\La$
	$$f_{i\nu(i)}(e_j) = m_{i\nu(i)}\delta_{ij}$$
By (c), this factors through $\Omega(W_0)$. 
Now define  $\phi(m_{i\nu(i)}) := f_{i\nu(i)}$. As in Lemma  \ref{lem:6.2}
	the map  $\phi$ is injective. 
The map 
	$\phi$ is surjective: By part (b), the set of all $f_{i\nu(i)}$ is a basis for ${\rm Hom}_{\La}(\Omega(W_0), M)$.
\end{proof}

\bigskip

\begin{proposition}\label{prop:6.4}  Assume $M$ is indecomposable and not projective. Assume $\ul{\rm Hom}(W, M)=0$. 
Then ${\rm top}(M)  \cong {\rm soc} (M)$.
\end{proposition}

\bigskip

\begin{proof}
The modules $W_0$ and $\Omega(W_0)$ are cyclic since the tops are multiplicity-free. Write $W_0 = \Theta \La$ and $\Omega(W_0) = \Psi\La$, here $\Theta$ and 
$\Psi$ are taken as elements in $\oplus_{i\in Q_0} e_i\La$. 

Since $\Omega^2(W_0)\cong W_0$ we have $\Theta\Psi=0= \Psi\Theta$, and there are 
exact sequences
$$0\to \Theta \La \to \La \to \Psi \La \to 0, \ \ \mbox{and} \ \ 
0\to \Psi \La \to \La \to \Theta \La \to 0.
$$
We apply the functor $(-, M):= {\rm Hom}_{\La}(-, M)$ to the first exact sequence, it takes it to an exact sequence
$$0\to (\Psi \La, M) \to (\La, M) \to (\Theta \La, M ) \to 0.$$
We identify the terms, as vector spaces. The middle is $M$. 
Furthermore
$$(\Psi \La, M) \cong \{ m\in M \mid m\Theta =0\}
\ \mbox{ and }
(\Theta\La, M) \cong \{ m\in M\mid m\Psi = 0\}
$$
where we view $\Theta$ and $\Psi$ as linear maps $M\to M$.
Hence we have an exact sequence
$$0 \to {\rm Ker}(\Theta) \to  M \to {\rm Ker}(\Psi) \to 0,
$$
which shows that  $M/{\rm Ker}(\Theta) \cong {\rm Ker}(\Psi)$.

\bigskip

By  Lemma \ref{lem:6.3}, ${\rm Ker}(\Theta) \cong {\rm rad} M$, and by Lemma \ref{lem:6.2}, we have 
${\rm Ker}(\Psi) \cong {\rm soc}(M)$. 
This shows that ${\rm top}(M) = M/{\rm rad}(M) \cong {\rm soc}(M)$ as
vector spaces, as required.
	\end{proof}

\bigskip

\begin{corollary} If $\ul{\rm Hom}(W, M)=0$ then 
	$M$ is $\Omega$-periodic.
	\end{corollary}

\begin{proof}  We have $\ul{\rm Hom}(W, M) \cong \ul{\rm Hom}(W, \Omega^n(M))$ for all $n\in \bZ$ since $W\cong \Omega(W)$. 
Hence by Proposition \ref{prop:6.4}, we have ${\rm top}(\Omega^n(M))\cong {\rm soc}(\Omega^n(M))$ for all $n\in \bZ$. 
Note that the top of $\Omega^r(M)$ is the socle of $\Omega^{r+1}(M)$. It  follows that the dimensions of the tops of the $\Omega^n(M)$ are constant and
therefore the dimensions of the $\Omega^n(M)$ are bounded.

Hence there is some integer $d$ such that infinitely many $\Omega^n(M)$ have
dimension $d$. Now we can apply \cite{CB} again which shows that some $\Omega^m(M)$ has 
$\tau$-period $1$, that is,  $\Omega$-period 2. Therefore $M$ is $\Omega$-periodic.
\end{proof}

\medskip

We conclude that 
on a component of a module $M$ which is not $\Omega$-periodic
the additive function $d_W$ above must be non-zero. Hence by Proposition \ref{prop:6.1} we can deduce 
the graph structure of a component.

\bigskip

\bigskip

\subsection{Auslander-Reiten components of simple modules and of some arrow modules}

In this part we assume that 
$H$ is a hybrid algebra (which may have virtual arrows), with distinguished set of triangles $\cT$, and we exclude the local algebra
with two virtual loops. We investigate
the position of simple modules, and of some modules generated by arrows, in the stable
AR quiver of $H$.
We say that a component is of type $A$ if its tree class is one of $A_{\infty}^{\infty}$ or $A_{\infty}$, or $\wt{A}_n$, or $A_n$ for some $n\geq 2$, and
we say it is of type  $D$ if its tree class
is one of $D_{\infty}$ or $\wt{D}_n$ or $D_n$.   For a vertex $i$ of $Q$ we denote the module ${\rm rad}(e_iH)/{\rm soc} (e_iH)$ by $M_i$ (the 'middle').
\medskip

\subsubsection{\bf Arrow modules for arrows not in $\cT$}

Take an arrow $\beta\not\in \cT$.   Then it is easy to see that
$\Omega^r(\beta H)\cong f^r(\beta)H$ for $r\geq 1$. Hence  $\beta H$ has $\Omega$-period equal to $r_{\beta}$ where
$r_{\beta}$ is the length of the $f$-orbit of $\beta$.
This is also true if some $f^s(\beta)$ is virtual of type (a), in which case the corresponding module is simple.  Furthermore all $\Omega$-translates are indecomposable and hence belong to ends of tubes
in the stable AR-quiver.

\medskip
\subsubsection{\bf Simple modules at biserial vertices, and at quaternion vertices}

(a) \ 
Assume $i$ is a biserial vertex. If there is a virtual loop at $i$ then by the previous, the simple module at $i$ is periodic at the end of a tube.
Now suppose the arrows starting at $i$ are not virtual. By Lemma \ref{lem:7.8} (see the Appendix), the 'middle' $M_i$ of $e_iH$ is the direct sum of two indecomposable modules.
Hence we have an almost split sequence $0\to \Omega(S_i)\to P_i \oplus M_i \to \Omega^{-1}(S_i)\to 0$ and $\Omega(S_i)$ has two predecessors in its stable component.
This could be in the middle of some component of type $A$  or possible in a component of type $D$  away from the edge.
In fact, it might even be in some tube. 

(b) \ If $i$ is a quaternion vertex,
with no singular relation close to $i$
(eg excluding  the singular tetrahedral, disc, triangle algebra)
then $S_i$ is periodic of period four. The proofs in \cite{WSA}, \cite{WSA-GV} and \cite{WSA-SD} generalize. This also works for \ref{ssec:4.2}(2c) and
for the algebra in \ref{ssec:4.6}.

\medskip

\subsubsection{\bf Simple modules at hybrid vertices}

\begin{lemma}\label{lem:6.6} Assume $H$ is a hybrid algebra but
	is not the algebra \ref{ssec:4.2}(2c) or the algebra \ref{ssec:4.6}.
Let $i$ be a vertex and $\alpha, \ba$ are arrows starting
at $i$ where $\alpha\in \cT$ and $\ba\not\in \cT$.
	Let $M:=  M_i =  {\rm rad}(e_iH)/{\rm soc}(e_iH)$.\\
(a) The module $M_i$ is indecomposable and it  occurs in two different AR-sequences as the non-projective middle term.\\
(b) If $H$ is not of finite type then the component of $S_i$ is of type $D$.
\end{lemma}

\begin{remark}\label{rem:6.2} Consider the algebra \ref{ssec:4.2}(2c), this
        has a hybrid vertex.
The algebra is special biserial (see Lemma \ref{lem:4.1}).
        Consider the simple module $S_2$ at the hybrid vertex, by Lemma \ref{lem:4.1} we know thar ${\rm rad}(e_2H)/S_2$ is the direct sum of two non-zero
        modules, and  $S_2$
        belongs to a component of tree class $A_{\infty}^{\infty}$.
        Similarly the algebra in \ref{ssec:4.6} has hybrid vertices $1, 3$
        but  the
        modules ${\rm rad}(e_iH)/S_i$ for $i=1, 3$ are
        decomposable. 
\end{remark}

\begin{proof} 
(a) \ Assume $\alpha \in \cT$ and $\ba \not\in \cT$. Note that then $f(\alpha)\neq \ba$.
As a preliminary part, we show that always $f(\alpha)f^2(\alpha)\ba =0$. 

	If now, then by (2) of Definition \ref{def:3.3} we have that
	$f(\alpha), g(\alpha) \in \cT$ and $g(\alpha)$ is virtual or critical.
	Suppose $g(\alpha)$ is virtual, then $n_{\alpha} = n_{g(\alpha)}\leq 2$.We cannot have $\alpha = g(\alpha)$ since this would imply $f(\alpha)= \ba$. 
So $g$ must have a 2-cycle $(\alpha \ g(\alpha))$, bu then $f(g(\alpha)) = \ba$.This gives a contradiction since with $g(\alpha)\in \cT$ also $f(g(\alpha))\in \cT$ but $\ba\not\in \cT$. 
	This shows that $g(\alpha)$ is not virtual.

Suppose $g(\alpha)$ is critical, consider first the case when the $g$-cycle of
	$g(\alpha)$ does not have a loop, then we use the diagram 3.1.1 with 
	$\tau = g(\alpha)$. Then $\xi$ must be virtual and therefore the arrow $y\to x$ must be in $\cT$, and then also $\ba$ is in $\cT$, a contradiction.
	Similarly one gets a contradiction in the other case, ie where
	$H$ is the algebra \ref{ssec:4.2}(2c). Hence $g(\alpha)$ is not critical.

\medskip

	The module  $M$ is indecomposable by Lemma \ref{lem:7.8}.
Therefore it is the indecomposable non-projective middle term of the AR-sequence starting with $\Omega(S_i)$.
Moreover we have a non-split short exact sequence
$$0\to V\to M\to U\to 0
\leqno{(*)}$$
where $V= \bar{\alpha} H/\langle B_{\ba}\rangle $ and $U = \alpha H/\langle A_{\ba}\rangle$.
Note that this is true also when $\ba$ is virtual.
We show first $V\cong \Omega^2(U)$, and next that ${\rm Ext}^1(U, V)\cong K$.
With these,  it will follow that (*) is an AR-sequence.
Let $j=t(\alpha)$ and $y= t(f(\alpha))$.

(i) We claim that $U$ is isomorphic to $e_jH/f(\alpha) H$:
Consider the projective cover $\pi: e_jH\to U$ given by $\pi(x) = \alpha x + \langle A_{\ba}\rangle$. Then
$\pi(f(\alpha))=0$ and hence $f(\alpha)H\subseteq {\rm Ker}(\pi)$.
	We can compare dimensions, applying Lemma \ref{lem:7.8}. 
The dimension of $U$ is $m_{\alpha}n_{\alpha}-1$ and
we have $\dim e_jH= m_{\alpha}n_{\alpha} + m_{f(\alpha)}n_{f(\alpha)}$. Hence
	the kernel of $\pi$ has dimension $n_{f(\alpha)}m_{f(\alpha)} + 1 = \dim f(\alpha) H$, and 
we have equality.
This implies that $\Omega(U) \cong f(\alpha)H$.

(ii) We claim that $\Omega(f(\alpha)H) \cong f^2(\alpha)\ba H$, and that it is isomorphic to $V$:
Let $\pi: e_yH\to f(\alpha)H$ be the projective cover, given by $e_yx\mapsto f(\alpha)x$.
As we have shown in the preliminary step, we always have
$f(\alpha)f^2(\alpha)\ba = 0$, so  
$f^2(\alpha)\ba H$ is contained in the kernel of $\pi$. By comparing dimensions we see that
it is equal.
To show that this is isomorphic to $V$, 
consider left multiplication with $f^2(\alpha)$ from $\ba H$ to $f^2(\alpha) \ba H$.
	This is a surjective  $H$-module homomorphism. By Lemma \ref{lem:7.8}, 
	$f^2(\alpha)\ba H$ has dimension
	$m_{f^2(\alpha)}n_{f^2(\alpha)}-1$ and 
	$\dim \ba H = m_{\ba} n_{\ba} = m_{f^2(\alpha)}n_{f^2(\alpha)}$
	noting $\ba = g(f^2(\alpha))$. So the kernel is  equal to $\langle B_{\ba}\rangle$.

(iii) \ It remains to show that ${\rm Ext}^1(U, V)$ is at most 1-dimensional (we know already that it is non-zero).
We have an exact sequence
$$ Ve_j\cong {\rm Hom}(e_jH, V)\stackrel{\iota^*}\to  {\rm Hom}(f(\alpha)H, V)\to {\rm Ext}^1(U, V)\to 0
$$
where $0 \to  f(\alpha)H\stackrel{\iota}\to e_jH$ is the inclusion map.

Assume first that $\ba$ is virtual. Then $V$ is 1-dimensional and spanned by the coset of $\alpha f(\alpha)$, so it is isomorphic to the simple
module $S_y$. In particular $Ve_y = V$ is 1-dimensional, and hence the quotient ${\rm Ext}^1(U, V)$ is at most 1-dimensional.

Now assume $\ba$ is not virtual. 
We have ${\rm Hom}(f(\alpha)H, V) \cong \{ v\in Ve_y \mid \ vf^2(\alpha)\ba  =0\}$.
The space $Ve_y$ is spanned by the  (cosets of) initial submonomials of $A_{\ba}$  which end at vertex $y$, that is which end in either
$f(\alpha)$ or in $\beta:= g^{-1}(f^2(\alpha))$.

Suppose $p$ is an initial submonomial of $A_{\ba}$ ending in $f(\alpha)$. By
the preliminary fact, we know that $p f^2(\alpha)\ba=0$, and  we deduce that there is a homomorphism  $\theta_p: f(\alpha)H \to V$ taking $f(\alpha)$ to $p$. We   claim that this
is in the image of $\iota^*$:
Such a monomial $p$ has a factorisation $p=\wt{p}\cdot f(\alpha)$ with $\wt{p}$ a monomial of positive length. There is a homomorphism $\wt{\theta}$
$e_jH\to V$ taking $e_j$ to $\wt{p}$ and hence
$\theta = \wt{\theta}\circ \iota$.

Now consider an initial  submonomial $p$ of $A_{\ba}$ ending in $\beta$.  If $p\neq A_{\ba}$ then 
$p f^2(\alpha)\ba$ is again an initial submonomial of $A_{\ba}$ and is non-zero in the algebra. This means that we do not have a homomorphism 
taking $f(\alpha)$ to $p$. 
This leaves only the case $p=A_{\alpha}$ so that the ext space is at most 1-dimensional. 
(In fact, this last case gives rise to the non-split short exact sequence).

(b) By assumption, $\ba H$ is $\Omega$-periodic. Let $W$ be the direct sum of
the distinct $\Omega$-translates of $\ba H$. Then $W\cong \Omega(W)$ and
$d_W(-)$ is an additive function on any non-periodic component
on which it does not vanish. By assumption, $H$ is of infinite type and then the summands of $W$ belong to tubes.
On the other hand, since
$H$ has infinite type, by part (a) the component of $M_i$ cannot
be a tube. The inclusion $\ba H \to \Omega(S_i)$ is nonzero in 
the stable category. Therefore  $d_W$ is non-zero on this component. 
We have $d_W(M_i) = 2d_W(\Omega(S_i)) \neq 0$ by exactness. Comparing with a
general additive function on components as described in \cite{HPR}, it follows that the component is of type $D$.
\end{proof}

\begin{remark} We see from the proof  $U$ or $V$ can be simple, or even both.
Consider the algebra
 $H$ with triangular quiver. We use the
	notation as in \ref{ssec:4.3}, and take  $\cT = \{ \alpha_i \}$ 
	and we take  $m_{\bullet}= (2, 1, 1)$. 
Then $\beta_2$ and $\beta_3$ are virtual, and we have
	$\Omega(S_1) \cong \Omega^{-1}(S_2)$. In this case, all three simple
	modules are of type $D$, in fact they are all in the same component
	which has tree class $\wt{D}_5$. Consider $M_3 = {\rm rad}(e_3H)/S_3$,
	in this case both $U$ and $V$ are simple.
\end{remark}

\vspace*{1cm}

\section{Appendix: Consistency, bases and dimensions}

This extends to the general case what was done for regular hybrid algebras in Section 2.
\medskip

\subsection{Consistency}

In this Section we assume throughout that $H$ is a hybrid algebra, which is not local, and is not an algebra considered in detail in Section 4.
With this assumption, we can use the diagrams in \ref{ssec:3.1}, see also
Corollary \ref{cor:3.9}.

\medskip

\begin{lemma} \label{lem:7.1} Assume $\ba$ is a virtual arrow, and $\alpha, \ba\in \cT$.  If $\ba$ is not a loop then
there are six  relations of type $\zeta$ or $\xi$  in which
$\ba$ occurs.
If $\ba$ is a loop then there are four relations of type $\zeta$ or $\xi$ in which $\ba$ occurs.
In both cases, each of these is zero in $H$.
\end{lemma}

%\medskip
%
%{\it Proof }
The proof is the same as that of Lemma 3.3 in   \cite{WSA-corr}, using the diagrams displayed in  \ref{ssec:3.1}. 
See also Corollary \ref{cor:3.9}.

\bigskip

\begin{lemma} \label{lem:7.2} 
 Assume $|A_{\alpha}|\geq 2$ but $\alpha$ is not critical. Let $\zeta = \zeta_{\alpha}:= \alpha f(\alpha)g(f(\alpha))$. \\
(a) If $\alpha, \ba\in \cT$ and $\ba$ is virtual or critical, then $\zeta\equiv A_{\alpha}$. Moreover
$$\zeta f^2(\ba) \equiv B_{\alpha}, \ \ \zeta g(f(\ba)) =0, \ \  g^{-1}(\alpha) \zeta  \equiv B_{g^{-1}(\alpha)}, \ \ f^{-1}(\alpha)\zeta = 0.
$$
Furthermore $B_{\alpha}J=0=JB_{\alpha}$ and $B_{g^{-1}(\alpha)}J=0=JB_{g^{-1}(\alpha)}$.\\
	(b) Otherwise $\zeta = 0$.
\end{lemma}

\bigskip

\begin{proof} Part (b) is a direct consequence 
	of part (2) in Definition \ref{def:3.3}.\\
(a) By the assumptions, $\alpha$ is not
	virtual or critical. We know from 
	 3.1.1 and 3.1.2  that
	$\zeta \equiv A_{\alpha}$.
It is clear that $\zeta f^2(\ba) \equiv B_{\alpha}$ and $g^{-1}(\alpha)\zeta \equiv B_{g^{-1}(\alpha)}$.
Furthermore,
since $\zeta \equiv A_{\alpha} \equiv \ba f(\ba)$, any monomial of length three having this as a factor, and
which has 'type $\zeta$ or type $\xi$' must be zero in $H$, by Lemma \ref{lem:7.1}. We will uses this throughout the
proof (without further comments).

(i) \  $\zeta g(f\ba)=0$ and $f^{-1}(\alpha)\zeta = 0$:
By the preamble,
$$f^{-1}(\alpha)\zeta =  g^{-1}(\ba)\zeta \equiv g^{-1}(\ba) \ba f(\ba) = 0 \ \ \mbox{ and } 
\zeta g(f(\ba)) \equiv  \ba f(\ba)g(f(\ba)) = 0.
$$
Note that these imply $B_{g^{-1}\alpha}\cdot g(f\ba)= 0$ and $f^{-1}(\alpha)B_{\alpha} = 0$.

(ii) \   $B_{\alpha}J=0 = JB_{\alpha}$:  first we have  
$$B_{\alpha}\ba \equiv A_{\alpha}g^{-1}(\alpha)\ba =  A'\beta g^{-1}(\alpha) \ba = A'\xi_{\beta} = 0
$$
where $\beta = g^{-2}(\alpha)$ that is $A'\beta = A_{\alpha}$ (which has length $\geq 2$ by assumption). 
Next, $B_{\alpha}\alpha \equiv B_{\ba}\alpha$. If $\ba$ is virtual we can write this as
$$B_{\ba}\alpha = \ba g(\ba)\alpha = \ba f^2(\alpha) \alpha = \xi_{\ba} = 0.$$
Suppose $\ba$ is critical, then we have, since $\xi_{\ba}=0$, 
$$B_{\ba}\alpha  = \ba g(\ba)g^2(\ba) \alpha = \ba g(\ba) f^2(\alpha)\alpha \equiv \ba g(\ba)A_{f(g\ba)} = 0.
$$
It remains to show $g^{-1}(\alpha) B_{\alpha} =0$ which is $\equiv g^{-1}(\alpha)B_{\ba}$. 
If $\ba$ is virtual we have
$$g^{-1}(\alpha)B_{\ba}  = g^{-1}(\alpha)\ba g(\ba) = \zeta_{g^{-1}(\alpha)} = 0$$
If $\ba$ is critical 
$$g^{-1}(\alpha)B_{\ba}  = g^{-1}(\alpha)\ba g(\ba)g^2(\ba) = \zeta_{g^{-1}(\alpha)} g^2(\ba) = 0$$
(iii) \  $B_{g^{-1}(\alpha)}J=0 = JB_{g^{-1}(\alpha)}$:  This is similar to (ii). We omit details.
\end{proof}

%\  First  $B_{g^{-1}(\alpha)}g^{-1}(\alpha) = g^{-1}(\alpha)B_{\alpha}=0$ by the previous part, and  
%$B_{g^{-1}(\alpha)}g(f(\ba))=0$ by part (i). 
%Next, 
%$$g^{-2}(\alpha)B_{g^{-1}(\alpha)} = g^{-2}(\alpha)g^{-1}(\alpha) A_{\alpha} \equiv g^{-2}(\alpha)g^{-1}(\alpha)\ba f(\ba) 
%= \xi_{g^{-2}(\alpha)} f(\ba) = 0.$$
%For the last step we must show $0 = f(\ba) B_{g^{-1}(\alpha)}= f(\ba) g^{-1}(\alpha) A_{\alpha}$.
%If $\ba$ is critical we write this as
%$$f(\ba)g^{-1}(\alpha) \alpha C = \zeta_{f(\ba)} C =0
%$$
%since $f(\ba)$ is virtual (here $C$ is a monomial of positive length).
% Finally, let $\ba$ be virtual, then also $g(\ba)$ is virtual and using this we have
%$f(\ba)g^{-1}(\alpha)A_{\alpha}  \equiv  g(\ba) A_{\alpha}  = g(\ba)\alpha g(\alpha) C = 0$ 
%(for a monomial $C$). 
%\end{proof}

\bigskip

\begin{lemma}\label{lem:7.3} Assume $\alpha$ is an arrow with $|A_{\alpha}|\geq 2$ but $\alpha$ not critical. Let $\xi =\xi_{\alpha}:= \alpha g(\alpha)f(g(\alpha))$. \\
(a) Suppose $\alpha, \ba\in \cT$ and $f(\alpha)$ is virtual. Then $\xi \equiv A_{\ba}$.
Moreover
$$g^{-1}(\alpha) \xi=0, \ \ f^2(\alpha)\xi = B_{f^2(\alpha)}, \ \  \xi f^2(\alpha) \ =  B_{\ba}, \ \  \xi f^2(g(\alpha))=0.
$$
We have $B_{\ba}J=0=JB_{\ba}$ and $B_{f^2(\alpha)}J=0=JB_{f^2(\alpha)}$.\\
(b) Suppose $\alpha, \ba\in \cT$ and $f(\alpha)$ is critical. Then $\xi \equiv A_{\alpha}$. Moreover
$$ \xi g^{-1}(\alpha) = B_{\alpha}, \ \ \xi g^{-1}(f(\alpha))=0, \ \ g^{-1}(\alpha)\xi  = B_{g^{-1}(\alpha)}, \ \ f^2(\alpha)\xi =0.
$$
We have $B_{\alpha}J=0=JB_{\alpha}$ and $B_{g^{-1}(\alpha)}J=0 = JB_{g^{-1}(\alpha)}$.\\
(c) Otherwise $\xi = 0$.
\end{lemma}

{\it Proof} \ This is similar to the proof of   Lemma \ref{lem:7.2}. We omit the details.
$\Box$

\medskip

The following deals with another special case.

\begin{lemma} \label{lem:7.4}
Assume that either  $\alpha$  is virtual and $\alpha\in \cT$, or
$\alpha$  is critical and $\alpha, g(\alpha)\in \cT$. Then $A_{\alpha}J = \langle B_{\alpha}, A_{\ba}\rangle$ and $A_{\alpha}J^2 = \langle B_{\alpha}\rangle$
and $B_{\alpha}J=0$.
\end{lemma}

\bigskip

\begin{proof}  Assume first that $\alpha$ is virtual, that is $\alpha = A_{\alpha}$ and $\ba\in \cT$. Then  $A_{\alpha}J = \langle \alpha g(\alpha), \ \alpha f(\alpha)\rangle \ = \ \langle B_{\alpha}, A_{\ba}\rangle$. 
By considering the diagrams in 3.1.2  we see that
 $\ba$ is not virtual or critical. We apply Lemma \ref{lem:7.2} with
$\alpha, \ba$ interchanged and get
$$A_{\ba} \equiv \zeta_{\ba}, \ \ A_{\ba}J = B_{\ba}, \ \ B_{\ba}J=0.$$
Therefore $A_{\alpha}J^2 = \langle A_{\ba}J \rangle = \langle B_{\ba}\rangle = \langle B_{\alpha}\rangle$. 

Now assume $\alpha$ is critical with $g(\alpha)\in \cT$. We have $A_{\alpha} = \alpha g(\alpha)$ and 
$$A_{\alpha}J = \langle \alpha g(\alpha)g^2(\alpha), \ \alpha g(\alpha)f(g(\alpha)) \rangle = \langle B_{\alpha}, \ \xi_{\alpha}\rangle,
$$
and we have  $\xi_{\alpha}\equiv A_{\ba}$ (see 3.1.2(1)(b)). By   Lemma \ref{lem:7.3}  we have that
$A_{\ba}J=B_{\ba}$, and $B_{\ba}J=0$ which implies the statement.
\end{proof}

\medskip

\begin{lemma} \label{lem:7.5} Assume $\alpha$ is any arrow, then\\
(i) \ $B_{\alpha}J=0$ and $JB_{\alpha}=0$. \\
(ii) \ $B_{\alpha}$ is non-zero.
\end{lemma}

\begin{proof} (i)  It suffices to show  that for an arbitrary arrow $\alpha$ we have $\alpha B_{g(\alpha)}$, that is $\alpha B_{f(\alpha)} = 0$. Then part (i)
follows using identity (3) of the definition \ref{def:3.3}, and 
	an identity such as $\alpha B_{g(\alpha)} = B_{\alpha}\alpha$.

If $\alpha\not\in \cT$ then $\alpha B_{f(\alpha)} =0$ by identity (1) of definition \ref{def:3.3}, so we assume now that
$\alpha \in \cT$. Then $f(\alpha)$ cannot be virtual of type (a) and therefore $|B_{f(\alpha)}|\geq 2$. 

(1) Assume $|B_{f(\alpha)}|=2$.  Then $\alpha B_{f(\alpha)} = \zeta_{\alpha}$. 
	If $f(\alpha)$ is virtual then $\zeta_{\alpha}=0$ by Lemma \ref{lem:7.1}. Assume now that $f(\alpha)$ is not virtual, it also is not critical (since
	$|B_{f(\alpha)}|\neq 3$). Therefore $\zeta_{\alpha}=0$ by 
	identity (2) of Definition \ref{def:3.3}.

(2) \ Assume $|B_{f(\alpha)}| = 3$, then 
$\alpha B_{f(\alpha)} = \zeta_{\alpha} g^{-1}(f(\alpha))$. This is zero 
unless $\ba\in \cT$ and $\ba$ is critical or virtual. Suppose $\ba$ is critical
	or virtual. Note first that we see from 3.1.1, 3.1.2 that $\alpha$ is not a loop. Therefore $\alpha B_{f(\alpha)}$ is not a cyclic path.
We also see from 3.1.1 and 3.1.2 that $\alpha$ cannot be virtual 
	or critical. That is, the assumption of Lemma \ref{lem:7.2} holds. 
It follows that  $\zeta_{\alpha}g^{-1}(f(\alpha))$ is zero by 7.2 
	(it is not cyclic and cannot be $\equiv B_{\alpha}$).

(3) Now assume $|B_{f(\alpha)}| \geq 4$. Then $\alpha B_{f(\alpha)} = \zeta_{\alpha}C$ where $C$ is a monomial of length $\geq 2$. Suppose
	$\ba$ is virtual or critical, then $\alpha$ is not virtual or critical 
	(see 3.1.1 or 3.1.2). By Lemma \ref{lem:7.2} we know $\zeta_{\alpha}J = \langle B_{\alpha}\rangle$ and $B_{\alpha}J =0$ and hence $\alpha B_{f(\alpha)}=0$.

	\bigskip

(ii) When the vertex $i=s(\alpha)$ is quaternion, the statement 
is proved in 4.5 of \cite{WSA-GV}.
Suppose $i$ is biserial. From the relations,
the only submonomials of $B_{\alpha}$ which occur in a minimal relation
are $B_{\alpha}$ itself and $A_{\alpha}$ and $A_{g(\alpha)}$.
In general, $A_{\alpha}$ occurs in a relation
$\ba f(\ba) - c_{\alpha}A_{\alpha}$ but this is not the case when
$i$ is biserial. Similarly
$A_{g(\alpha)}$ could occur in a relation $f(\alpha)f^2(\alpha) - c_{g(\alpha)}A_{g(\alpha)}$ but not if $i$ is biserial since in that case
$f(\alpha)f^2(\alpha)$ is zero (or a scalar multiple of $B_{f(\alpha)}$).
Hence $B_{\alpha}$ is non-zero in $H$.

Now assume that $i$ is hybrid, say $\alpha  \in \cT$ and $\ba \not\in \cT$.
Then $A_{\alpha}$ does not occur in a defining relation. We have the relation
$f(\alpha)f^2(\alpha) = c_{\alpha}A_{g(\alpha)}$ but this does not give
a relation which forces $B_{\alpha}$ to be zero in $H$.
\end{proof}

\medskip

\begin{lemma}\label{lem:7.6}
Consider a path of length four  of the form
$p:= \alpha g(\alpha)\beta g(\beta)$ where $\beta= f(g(\alpha))$.\\
(a). If  $f(\alpha)$ is virtual or critical
then $p$ is a non-zero scalar multiple of $B_{\alpha}$. \\
(b) \ Otherwise it is zero.
\end{lemma}

\begin{proof}
We can write $p=\alpha \zeta_{g(\alpha)}$ and also $p=\xi_{\alpha} g(\beta)$.
By Lemma \ref{lem:7.2} we know $\alpha \zeta_{g(\alpha)}\neq 0$ if and only if $f(\alpha)$ is virtual or
critical and if so it is $\equiv B_{\alpha}$ (which is $\equiv B_{\ba}$).

By Lemma \ref{lem:7.3} we have $\xi_{\alpha} g(\beta)\neq 0$ if and only if $f(\alpha)$ is virtual or critical, and if so then $p$ is
a cyclic path of length four  $\equiv B_{\alpha}$.
 \end{proof}

\bigskip

\subsection{Bases and dimension}

In the following write $|A_{\alpha}|=\ell$ and
$|A_{\ba}|=\bar{\ell}$. We also write $[A_{\alpha}]_j$ for the initial submonomial of $A_{\alpha}$ of length $j$.

\bigskip

\begin{lemma}\label{lem:7.7} Assume $\alpha$ is an arrow of $Q$. Then the set
$\{ [A_{\alpha}]_j \mid 1\leq j\leq \ell,  \ B_{\alpha} \}$ is linearly independent.
 \end{lemma}

\bigskip
This is similar to proofs in Section 2, we omit details.

\medskip

\begin{lemma}\label{lem:7.8}
Assume $i$ is a vertex which is either  biserial or hybrid. Then \\
(a) \ $e_iH$ has basis consisting of all proper initial submonomials of $B_{\alpha}, B_{\ba}$ together with $e_i$ and $B_{\alpha}$.\\
(b) \ $\dim e_iH = m_{\alpha}n_{\alpha} + m_{\ba}n_{\ba}$.\\
	(c) \ If $\alpha \in \cT$ then $\dim \alpha H = m_{\alpha}n_{\alpha}+1$,
	otherwise $\alpha H$ has dimension $m_{\alpha}n_{\alpha}$. The module
	$\alpha g(\alpha)H$ has dimension $m_{\alpha}n_{\alpha}-1$.\\
	(d) \ Let $M_i= {\rm rad}(e_iH)/{\rm soc}(e_iH)$. If $i$ is biserial and $\alpha, \ba$ are not virtual of type (a) then  $M_i$ is the direct sum of
two uniserial modules. If one of $\alpha, \ba$ is virtual of type (a) then $M_i$ is uniserial.
 If $i$ is hybrid, and
	$H$ is not the algebra in \ref{ssec:4.2}(2c) or \ref{ssec:4.6} then
	$M_i$ is indecomposable.
\end{lemma}

\bigskip

\begin{proof}  We prove part (a), then parts (b) and (c) follow directly.
We may assume $\ba \not\in \cT$. 
The
given set spans $e_iH$ by Lemmas \ref{lem:7.2} and \ref{lem:7.3}.
We show  linear independence.
Take a linear combination
$$\sum_{j=1}^{\ell} a_j[A_{\alpha}]_j + \sum_{t=1}^{\bar{\ell}} d_t [A_{\ba}]_t + sB_{\alpha} = 0. \leqno{(*)}
$$
Let $\beta = f^{-1}(\ba)=g^{-1}(\alpha)$,  then $\beta\ba =0$, unless possibly $\beta = f(\ba)$, a loop, and
$\beta\ba = b_{\beta}B_{\beta}$. But then noting that $\alpha\in \cT$ it follows that $|Q_0|\leq 2$, which we have excluded.
Therefore we have $\beta\ba = 0$. We premultiply (*) with $\beta$ and obtain
$\sum_{j=1}^{\ell} a_j[A_{\beta}]_{j+1} = 0$, 
and by Lemma \ref{lem:7.7} it follows that $a_j=0$ for $1\leq j\leq \ell$.
Now applying Lemma \ref{lem:7.7} again implies $d_t=0$ for all $t$, and $s=0$.

(d) When $i$ is biserial, the claim also follows from part (a). Now suppose $i$ is a hybrid vertex, so $\alpha\not\in \cT$.
If $f(\alpha)$ is not virtual then $M_i$ can be viewed as a string module (see \cite[II.3]{E1}),  hence it is indecomposable. 
If $f(\alpha)$ is virtual then by the assumption that $H$ is not the algebra in \ref{ssec:4.2}(2c) or \ref{ssec:4.6} one checks that  at least one of
$A_{\alpha}$ and  $A_{\ba}$ has length  $> 2$, and then one verifies directly that $M_i$ is indecomposable.
\end{proof}

\bigskip

\begin{lemma}\label{lem:7.9}
 Assume $i$ is a periodic vertex so that  $\alpha, \ba$ are both in $\cT$.
 Then the set $\{ [A_{\alpha}]_j, (j\leq \ell),  [A_{\ba}]_t,  (t\leq \bar{\ell}),  B_{\alpha}\}$ is linearly independent,
 except when $H$ is the singular spherical algebra, or $H$ is the singular triangle algebra.
 \end{lemma}

\medskip

 This is proved in \cite{WSA-GV} (see Proposition 4.9).

\end{document}